\documentclass[11 pt]{amsart}
 \usepackage{amsmath}
\usepackage{amssymb}
\usepackage{amsfonts}

 \usepackage{eucal}
      \makeatletter
     \def\section{\@startsection{section}{1}%
     \z@{.7\linespacing\@plus\linespacing}{.5\linespacing}%
     {\bfseries
     \centering
     }}
     \def\@secnumfont{\bfseries}
     \makeatother
  \usepackage{a4wide}
\newtheorem{theorem}{Theorem}[section]
\newtheorem{lemma}[theorem]{Lemma}
\newtheorem{corollary}[theorem]{Corollary}
\newtheorem{proposition}[theorem]{Proposition}

\theoremstyle{definition}

\newtheorem{example}[theorem]{Example}

\theoremstyle{remarks}
\newtheorem{remarks}[theorem]{\bf Remarks}
\theoremstyle{remark}
\newtheorem{remark}[theorem]{\bf Remark}

\theoremstyle{problem}
\newtheorem{problem}[theorem]{\bf Problem}

\numberwithin{equation}{section}

\def\ddate {\sevenrm \ifcase\month\or January\or
February\or March\or April\or May\or June\or July\or
August\or September\or October\or November\or December\fi\! {\the\day}, \!{\sevenrm\the\year}}

\def \a{{\alpha}}
\def \b{{\beta}}
\def \D{{\Delta}}
\def \d{{\delta}}
\def \e{{\varepsilon}}
\def \g{{\gamma}}

\def \k{{\kappa}}
\def \l{{\lambda}}

\def \O{{\Omega}}
\def \p{{\varphi}}
\def \t{{\vartheta}}

\def \m{{\mu}}
\def \s{{\sigma}}

\def \noi{{\noindent}}
\def \qq{{\qquad}}
\def \dd{{\rm d}}
\def\F{{\mathcal F}}  
  
\def\E{{\mathbb E \,}}

\def\R{{\mathbb R}}

\def\Z{{\mathbb Z}}
 
\def\Q{{\mathbb Q}}

\def\N{{\mathbb N}}

\def\C{{\mathbb C}}  
  
  \scrollmode
 \font\sevenrm= cmr10 at 7 pt

\begin{document} 
\title[Farey  sequence and  
 quadratic Farey sums]
{\bf On 
the Uniform Distribution  (mod 1) of the Farey  Sequence,   quadratic Farey  and Riemann sums
with a remark  on  local integrals of $\boldsymbol{\zeta(s)}$
}

\author{Michel J.\,G. Weber}

\address{IRMA, 10 rue du G\'en\'eral Zimmer,
67084 Strasbourg Cedex, France}
\email{michel.weber@math.unistra.fr, m.j.g.weber@gmail.com}
\renewcommand{\thefootnote}{} {{
\footnote{2010 \emph{Mathematics Subject Classification}. Primary 42A75,    Secondary 42A24, 42B25.}
\footnote{\emph{Key words and phrases}: Farey fractions, Farey sums, Riemann sums, Riemann Hypothesis, uniform distribution,  Riemann zeta function, local integrals, Stepanov's norm.}
\setcounter{footnote}{0}

\maketitle 

 \begin{abstract} For   $1$-periodic functions $f$   satisfying only a weak local regularity assumption of Dini's type at  rational points of $]0,1[$, we study  the Farey sums
 $$F_n(f)= \sum_{\frac{\k}{\l}\in \F_n} f\big(\frac{\k}{\l}\big),\qq F_{n,\s}(f)= \sum_{\frac{\k}{\l}\in \F_n} \frac{1}{\k^\s\l^\s}f\big(\frac{\k}{\l}\big),\qq  1/2\le \s<1  ,
$$
where $\F_n$
    is the Farey series of order $n\ge 1$. We obtain  sharp estimates of $F_{n,\s}(f)$, for  all $0< \s\le1$. We prove similar results for  the corresponding Riemann quadratic sums
$$ S_{n,\s}(f) \ =\  \sum_{1\le k\le \ell \le n}\frac{1}{(k\ell)^{\s }}\, f\big( \frac{k}{\ell}\big). 
$$
These sums are related to  local integrals of the  Riemann  zeta-function 
 over bounded intervals $I$, which are considered in the last part of the paper.
\end{abstract}
 


\section{\bf Introduction}\label{s0}
 
 Let
    $\F_n=\big\{ \frac{j}{m}, 1\le j\le m\le n, (j,m)=1\big\}$
    be the Farey series of order $n\ge 1$.
Let also $f:\Q\cap [0,1]\to \C$  arbitrary and $\s\ge 0$. In this work we study the  Farey sums   
\begin{equation} \label{F_n(f)}F_n(f)= \sum_{\frac{\k}{\l}\in \F_n} f\big(\frac{\k}{\l}\big), \qq F_{n,\s}(f)= \sum_{\frac{\k}{\l}\in \F_n} \frac{1}{\k^\s\l^\s}f\big(\frac{\k}{\l}\big)\qq n=1,2,\ldots.
\end{equation}
As $F_{n,0}(f)=F_n(f)$,
the second sums    generalize the first ones. 
\vskip 3 pt   We also study  the corresponding Riemann quadratic sums
\begin{eqnarray}\label{S_n(f)} S_{n,\s}(f) \ =\  \sum_{1\le k\le \ell \le n}\frac{1}{(k\ell)^{\s }}\, f\Big( \frac{k}{\ell}\Big) 
\end{eqnarray}
where $0\le \s <1$ and $f\in L^0([0,1])$. These sums  have a simpler structure as being  weightings  of the  Riemann sums of 
$f_\s(x)= 
f(x)/ x^{\s}$, and are connected with $\zeta$-local integrals, where  
$\zeta$ is the Riemann-zeta function, since for instance (see section \ref{s4})
\begin{equation*}\int_n^{n+1} |\zeta(\s +it ) |^2 \dd t\,=\, \hbox{{\it Quadratic Riemann sum} + {\it Error term}}\,.
\end{equation*}
   \vskip 3 pt Some preliminary considerations are necessary.
 It is well-known that the Farey fractions are uniformly distributed (mod 1), see   Mikol\'as \cite{M1}, whence by the Weyl criterion,  for any Riemann integrable function on $[0,1]$,
\begin{equation}\label{limF} F_{n}(f)\sim \Phi(n)\, \int_0^1 f(x)\dd x ,\qq \quad n\to \infty.\end{equation}
Here $\Phi(n)= \#(\F_n)$ and  $\Phi(n)=\sum_{d=1}^n \p (d) 
\sim \frac{3}{\pi^2} n^2$,     $ \p (d)$ being Euler's totient function. The problem of estimating the error term
\begin{equation}\label{Enf} E_n(f)\ =\ F_{n}(f)-\Phi(n)\,\int_0^1 f(x)\dd x ,
\end{equation}
 is connected with the Riemann Hypothesis, and was  studied by Mikol\'as \cite{M1} and by several authors, notably Codec\`a and Perelli \cite{CP},   see  references therein, and Yoshimoto \cite{Y}.  
  Farey sums and Riemann sums are linked by the formula   (\cite[Lemma\,2]{M1}),
\begin{equation}\label{farey.riemann}F_n(f)=\sum_{d\le n}D_{f}(d)\,M\big(\frac{n}{d}\big)=\sum_{d\le n} \m*D_f(d), \qq D_{f}(\ell)\, =\, \sum_{1\le k \le \ell}  f\Big(  \frac{k}{\ell}\Big),
\end{equation} 
 where  $*$ stands for the Dirichlet convolution product. 
One notes that $\Phi(n) =\sum_{d=1}^n dM(\frac{n}{d})$ where  $M(x)= \sum_{ \l \le x}\m(\l)$,  $\m(n)$  being the M\"obius function, and that $M(n)$ is the Farey sum
$$M(n)= \sum_{\frac{\k}{\l}\in \F_n} \cos 2\pi(\frac{\k}{\l}),$$
 which easily 
 follows
  from \eqref{inversemu}. By a result of  Littlewood \cite{L},  the Riemann Hypothesis is equivalent to
the assertion  \begin{equation}\label{RHM}M(x)=\mathcal O_\e\big(x^{\frac12 +\e}\big). \end{equation}
The  simplest example of a smooth periodic function $f(x) = \cos 2\pi x$,  thus  shows that the problem of estimating $E_{n}(f)$ ($=F_n(f)$ here)  is out of reach,   advances in this domain are therefore difficult.  
Farey sums much differ at this regard from Riemann sums  $R_{f}(\ell)= \frac{1}{\ell} D_{f}(\ell)$, 
since by a result of Wintner \cite[\S\,12]{Wi},  a continuous $1$-periodic function is analytic if and only if  there exists $q$, $0<q<1$, such that 
\begin{eqnarray}\label{Riemann.Wintner} R_{f}(\ell)\, =\,\int_0^1 f(t)\dd t  + \mathcal O( q^\ell).
\end{eqnarray}
A rate of convergence can be assigned, and the convergence of Riemann sums turns up the more rapid, the smoother $f$ is.
If $f$ is only Lebesgue integrable, the corresponding   convergence problems 
of   Riemann sums, and by extension Riemann equidistant sums,  Farey sums,  are  another attracting and  difficult  matter. We refer to Ch.\,XI of our book \cite{W1}.
\vskip 2 pt For the case considered above ($h(x)=\cos 2\pi x$, $F_{n,\s}=F_{n,\s}(h)$), we will prove that 
   \begin{equation}\label{Ih}
 F_{n,\s}\,=\, \frac{ n^{2(1-\s)} } {2(1-\s)\zeta(2)}\int_{0}^{1} \frac{\cos 2\pi t}{t^\s}\dd t\, +  \mathcal O_\s\big(n^{1-\s}\big)\qq \quad (0< \s<1).
\end{equation}
\vskip 3 pt 
The analogous formula to \eqref{farey.riemann} for the Farey sums $F_{n,\s}(f) $
is
\begin{equation}\label{fnf}
  F_{n,\s}(f)  \ =\ \sum_{d =1}^n\frac{1}{ d^{2\s}}\,\m*D_{f_\s}(d)\ =\ \sum_{d =1}^n \frac{D_{f_\s}(d)}{d^{2\s}}
 \, \Big(\,\sum_{\l\le \frac{n}{d}} \, \frac{\m( \l)} { \l^{2\s}}\,\Big).
\end{equation}
See  Lemma \ref{l2}. 
In comparison with \eqref{RHM} one knows (Mikol\'as \cite[Lemma\,8]{M1}) that for $\tau \ge 1/2$, 
RH is also equivalent to \begin{equation}\label{wmf}
 \sum_{\l\le n} \, \frac{\m( \l)} { \l^{\tau}}=\mathcal O_\e\big(n^{\frac12 +\e}\big)  .
\end{equation}   Further the Dirichlet series associated with $\m*D_{f_\s}$ being 
the product 
$$\frac{1}{\zeta(s)} \ \sum_{k=1}^\infty \frac{D_{f_\s}(k)}{k^s}.$$
  $ F_{n,\s}(f)$ can be precisely estimated  by using  Perron's formula, once estimates of $D_{f_\s}(k)$ are at disposal. This was used in Mikol\'as \cite{M1}  and   Codec\`a-Perelli \cite{CP}.
\vskip 3pt

The first formula in \eqref{farey.riemann} together with  \eqref{RHM} imply  that 
 \begin{equation}\label{farey.riemann1}F_n(f)=\mathcal O_\e\Big(n^{\frac12 +\e}\sum_{d\le n}\frac{|D_{f}(d)|}{d^{\frac12 +\e}}\Big),
\end{equation} 
for all $f\colon[0,1]\to \C$, assuming the validity of the RH. 
\vskip 2 pt Conversely we prove the following unconditional result, of very close order of magnitude.
\begin{theorem}\label{fnfomega} 
For
for infinitely many $n$, there exists  $f\colon [0,1]\to \C$ such that \begin{equation}\label{farey.riemann2}F_n(f)\,\ge\, n^{\frac12}\, \sum_{d\le n}\frac{|D_{f}(d)|}{d^{\frac12 }}.
\end{equation}
\end{theorem}
The proof is a combination of a theorem  of Pintz \cite{Pi}, which in particular implies that
\begin{equation} \label{omegamu}M(x) = \O\,\big(x^{\frac12}\big),
\end{equation}
and  infinite M\"obius inversion formula. \vskip 5 pt 
The work made in \cite{CP} concerns absolutely continuous functions   on $[0,1]$, or equivalently, continuous  functions $f$ with a derivative $f'$ almost everywhere and  $f'$  is Lebesgue integrable. It is further imposed that  $f'\in L^p[0,1]$ for some $p\in(1,2]$.  Let   $\mathcal C_{\rm cp}$ denotes this class of   functions. The main results obtained being of conditional nature, are  by definition   ineffective.
 However these results are nearly optimal with respect to the class $\mathcal C_{\rm cp}$, and \cite{CP} is one of the  central papers in the theory with Mikol\'as \cite{M1,M2},  notably by the ideas implemented. Some of these results were slightly extended in Yoshimoto \cite{Y}, who notably much investigated some specific remarkable  classes of functions related to Riemann sums.
  Let $\zeta(s)$ denotes  the Riemann zeta function. More precisely, it is assumed that  a weaker form of Riemann Hypothesis, noted RH($\a$) and  meaning that: $\s>\a \  \Rightarrow \   \zeta (\s+it)\neq 0$,
where $\a\in[1/2,1)$, 
  holds true.   One notes that then  $M(x)=\mathcal O_\e(x^{\a+\e})$. For $f\in \mathcal C_{\rm cp}$, there is a close link between the deviation of $R_f(n)$ from $\int_0^1 f(t)\dd t$, and Fourier coefficients of $f'$, which  is the basis of the approach 
 in \cite{CP}. 
 
 \begin{remark}\label{Kubert} 
  The following class of continuous functions $f\colon(0,1)\to \C$ satisfying Kubert identities
 $ f(x)= m^{s-1}\sum_{k=0}^{m-1}f\big( \frac{x+k}{m}\big) $ 
  for every $m\in \N$ and every $x\in(0,1)$, was studied in  \cite{Y}, and more recently in other papers. We study this class  of functions in a separate work, as well as some variant  in Carlitz \cite{Ca}.
 \end{remark}
  \vskip 3 pt 

In this paper we are interested in the study of Farey sums \eqref{F_n(f)} under minimal conditions,  and our results are effective and will  depend on the Fourier coefficients of $f$. Recall   that every function defined almost everywhere in $[0,1]$, in particular every integrable function,  has its Fourier series, see  Zygmund \cite[p.\,9]{Z}. 
 One motivation is that the limit \eqref{limF}  actually holds true for any function $f$ defined on rational points of $[0,1]$, whose  Riemann sums are converging, see \cite[Th.\,2]{M1}.
Another motivation lies in the fact that there are  important classes of functions $f$ having a (non integrable) derivative on $]0,1[$,  so $f$ is  not  absolutely continuous.
   The following   specific functions, familiar in Fourier analysis and relevant  in section \ref{s4}, are typical  cases,
 \begin{equation}\label{typicalcases}    \hbox{$g(a,x) =\frac{\sin a\log x}{  \log x}$,  $0<x<1$, \quad $g(a,1) =a$, $g(a,0)=0$, \quad ($a\neq0$) 
}.
 \end{equation} 

They  are not absolutely continuous on $(0,1]$, for otherwise they would have an absolutely continuous extension on $[0,1]$, thereby continuous on $[0,1]$, which is not. 
 \vskip 2 pt Let further  $g_\s(a,x)=g(a,x)/x^\s$, $0<\s<1$. Then  $g_\s(a,.)\in L^p[0,1]$  if $p <1/\s$,  $g_\s(a,.)'$ is not integrable, and as $g_\s(a,.)$ is unbounded in the neighbourhood of $0 $,   its  Fourier series does not converge absolutely.  
 
  For the problem studied, considering the restriction of  $g_\s(a,.)$ on $[\e,1]$, $\e>0$ small, and applying Euler-McLaurin's formula + Parseval's formula   is inoperant.
   
   We  will use the simplified notation (section \ref{s3})
   \begin{eqnarray}\label{g}   g(x)\ =\  g(1,x).
\end{eqnarray}

\vskip 7pt
 
 We   consider  in this work  $1$-periodic functions $f$ satisfying a weak local regularity assumption of Dini's type at   rational points of $]0,1[$, which is  in accordance with  the fact that Farey sums $F_n(f)$, $F_{n,\s}(f)$ are determined by the values taken by $f$ on rational numbers. More precisely,  introduce the class $\mathcal C$ of functions $f:[0,1]\to \R$   such that
\begin{align} \label{hyp.dini}     \hbox{\it   $f(x-0)=f(x+0)$ and $\int_{0+} \frac{| f(x+u)+f(x-u)-2 f(x)| }{u}\, \dd u<\infty$,  for each $x\in\Q\cap ]0,1[$.} 
\end{align}

This defines a fairly wide setting, one has  the following obvious inclusions: $ \mathcal C_1 \subset  \mathcal C_2\subset  \mathcal C$, where 
 \begin{eqnarray} \label{Class}& &  \mathcal C_1=\big\{f:[0,1]\to \R: \hbox{$f(x)$  is derivable for each  $x\in\Q\cap ]0,1[$\big\}}
  \cr & &  \mathcal C_2=\big\{f:[0,1]\to \R: \hbox{${f}(x)$ satisfies a Lipschitz condition of order $\a>0$},
 \cr & & \qq \qq \qq \hbox{in the neighbourhood of each $x\in\Q\cap ]0,1[$}\big\}.
 \end{eqnarray}

These functions being not necessarily absolutely continuous,   Euler-McLaurin's formula, which is the pivot   of the approach   in \cite{CP}, does not apply.  As      $x=0$ is excluded in definition \eqref{hyp.dini},  we note that $g_\s(a,.)\in \mathcal C$.
 \vskip 7 pt  
We clarify that  the approach used and most of the results  obtained extend with no difficulty to   classes of functions subject to 
sharper types of criteria such as the one of Jordan, Young, de la Vall\'ee-Poussin, Lebesgue, see  \cite[Vol.\,I]{B}. These ones being more elaborated we chosed to develop the present  work in this  simpler setting.

\vskip 7 pt The paper is organized as follows. In the two next sections we respectively state and give the proofs of results concerning Farey sums $F_{n}(f)$, $F_{n,\s}(f)$, and further comment and discuss our assumptions, comparing them notably with those in \cite{CP}. Section \ref{s1} contains preparatory  results which are interesting on their own. In section \ref{proof{t1a}}, the proofs of Theorems \ref{fnfomega} and  \ref{t1a} are given. Our results concerning quadratic Riemann sums $S_{n,\s}(f)$ are stated and proved in section \ref{s3}. In section \ref{s4},
we discuss some questions related to the  previous sections and concerning  local integrals of $\zeta(\s+ it)$.

The investigation of a related non trivial question concerning  the unboundedness of the $\mathcal S^2$-Stepanov's norm of  the Riemann zeta function $\zeta(\s+ it)$ for $1/2<\s<1$ is concluding the paper, the case $\s=1/2$ being trivial.


  \section{Farey sums and quadratic Farey  sums.}\label{s2}
 
 We first  provide general explicit formulas of  $F_{n,\s}(f)$ or $F_{n}(f)$, valid for all $f$ such that ${f_\s}\in \mathcal C$, with no additional condition. 
We next study these sums, mainly under two type of conditions on the complex Fourier coefficients of $f_\s$ where $f_\s(x)=\frac{f(x)}{x^\s}$, $\s \ge 0$.  We assume that: {\it   Either {\rm (i)}  the series  $\sum_{\ell \in \Z} c_{f_\s}(\ell)$ converges, or {\rm (ii)} the series $\sum_{\ell \in \Z} |c_{f_\s}(\ell)|$, converges.}
\vskip 3 pt
Recall some classical facts.  The first assumption  of course holds  if (for instance) ${f_\s}(x) $ is derivable at $x=0$.
   By Bernstein's theorem \cite[Vol.\,I,\,p.\,216]{B}, the second assumption holds   if $f\in {\rm Lip}(\a)$ for $\a>1/2$. It  is also implied by  the absolute convergence of the Fourier series
\begin{equation}\label{sfx}
\s_{f}(x)=\sum_{\nu\in \Z} c_{f}(\nu)e_\nu(x),\qq \quad (e_\nu(x)=e^{2{\rm i}\pi\nu x})
\end{equation}
at a single point, and   implies the absolute and uniform convergence of the series $\s_{f}(x)$ for all $x$. It also implies that $\s_{f}(x)$  converges to $f(x)$ at almost  every $x$,  and if $f$ is continuous on $[0,1]$, that convergence holds for every $x$. Further, if $f$ is absolutely continuous and $|f'|\ln^+ |f'|\in L^1([0,1]$, then $\s(f)$ converges absolutely. Also there exists an absolutely continuous function $f $ such that $\s(f)$ does not possess a single point of absolute convergence. See \cite[Vol.\,II,\,p.\,162]{B}.
\vskip 2 pt  Furthermore the second assumption  implies  that $f$ is bounded, since  using for instance  Riesz's criterion \cite[Vol.\,II]{B}, p.\,184, as the series $\s_{f}$ converges absolutely,  $f$ can be represented in the form
 $$f(x)= \int_{0}^{1}u(2\pi(x-t))v(2\pi t) \dd t ,$$ 
  $u,v\in L^2[0,1]$.
  By applying Young's inequality
  $$ \|u*v\|_r\le \|u\|_p\,\|v\|_q, \qq 1\le p,q, r\le \infty,\quad \frac{1}{r}=\frac{1}{p}+\frac{1}{q}-1,$$
   with $p=q=2$, $r=\infty$,   we deduce that $\|f\|_\infty <\infty$. 
    \vskip 3 pt
  In particular for $g_\s$ defined in \eqref{g}, we note that
  \begin{equation}\label{sfsigma}\sum_{\ell \in \Z} |c_{g_\s}(\ell)|=\infty  .
  \end{equation} 
\vskip 5 pt Consider first the case when  $f$ is such that ${f_\s}\in \mathcal C$, and no additional condition is imposed. 
    \begin{proposition}\label{Fp1} {\rm (i)} Let $f:[0,1]\to \R$ and assume that ${f_\s}\in \mathcal C$. 
Then we have,
\begin{align*} 
F_{n,\s}(f)
=\ & {f_\s}(1)\Big(\sum_{ d\le n}\frac{1}{d^{2\s}} \big(\sum_{\l\le \frac{n}{d}} \frac{\m(\l)}{\l^{2\s}}\big)\Big)
+c_{f_\s}(0)\Big[\sum_{ d\le n}\frac{d-1}{d^{2\s}} \Big(\sum_{\l\le \frac{n}{d}} \frac{\m(\l)}{\l^{2\s}}\Big)
\Big]
   \cr &\qq+
   \sum_{\ell\in \Z_*} c_{f_\s}(\ell) \Big[\sum_{ d\le n}\frac{1}{d^{2\s}}    \e_\ell(d)  \Big(\sum_{\l\le \frac{n}{d}} \frac{\m(\l)}{\l^{2\s}}\Big)\Big],
\end{align*}
where  
\begin{equation}\label{eps}\e_\ell (d)=\begin{cases} d-1 &\quad {\rm if}\ d|\ell, \cr
-1&\quad {\rm if}\ d\not|\ell.\end{cases}
\end{equation} 

{\rm (ii)}   Let $f:[0,1]\to \R$ be such that ${f }\in \mathcal C$.
Then we have,
 \begin{eqnarray*} 
E_n(f)\, =\, f(1)-\int_0^1f(x)\dd x 
+
   \sum_{\ell\in \Z_*} c_f(\ell)  \Big[\sum_{ d\le n}\e_\ell(d)  M \big(\frac{n}{d} \big)\Big].
\end{eqnarray*}
\end{proposition}

In the next theorems, we derive precise estimates of $F_{n,\s}(f)$  under the afore mentionned assumptions.

\begin{theorem}[$0<\s<1 $]\label{fp1}    Let $f:[0,1]\to \R$ be such that ${f_\s}\in \mathcal C$.
\vskip 3 pt   {\rm (i)} Assume that the series $\sum_{\ell \in \Z} c_{f_\s}(\ell)$ is convergent. Then
\begin{eqnarray*}  
   F_{n,\s}(f)& =&    \Big( \frac{\int_0^1f_\s(x) \dd x} {2(1-\s)\zeta(2)}  \Big)   n^{2(1-\s)}+ \D_{n,\s}(f)   ,
\end{eqnarray*}
     \begin{equation*}
\hbox{where}\ \ \D_{n,\s}(  f)\, =\,  \mathcal O (n^{1-2\s}\log n   )+ A (f_\s)+  \sum_{\ell\in \Z_* } c_{f_\s}(\ell)\Big(\sum_{ d\le n\atop d|\ell}\frac{1}{d^{2\s-1}}\sum_{\l\le \frac{n}{d}} \frac{\m(\l)}{\l^{2\s}}\Big),  
\end{equation*}
and  the constant $A (f_\s)$  is defined as follows,  
\begin{equation*} A (f_\s)=\begin{cases}{\ f_\s}(1)+    \frac{\zeta(2\s-1)}{\zeta(2\s)} \int_0^1f_\s(x)\dd x  -\displaystyle{\sum_{\ell\in \Z }} c_{f_\s}(\ell)\ & \ \hbox{if \ $\frac{1}{2}<\s<1$,} 
\cr \  0\ & \ \hbox{if \ $\s =\frac{1}{2}$,} 
\cr     -\displaystyle{\sum_{\ell\in \Z_* }}c_{f_\s}(\ell)\ & \ \hbox{if\ \ $0<\s<\frac{1}{2}$.} 
\end{cases}
\end{equation*}

\vskip 3 pt   {\rm (ii)}   Assume that the series $\sum_{\ell\in \Z_* } |c_{f_\s}(\ell)| $ is convergent. Then
\begin{eqnarray*}  
   F_{n,\s}(f)& =&    \Big( \frac{\int_0^1f_\s(x) \dd x} {2(1-\s)\zeta(2)}  \Big)   n^{2(1-\s)}+ \D_{n,\s}(f)   ,
\end{eqnarray*}
     where
      \begin{equation*}
 \D_{n,\s}(  f)\, =\, \begin{cases}
 A (f_\s)+ \mathcal O\Big(n^{1-2\s}\log n +  \displaystyle{\sum_{\ell\in \Z_* }} |c_{f_\s}(\ell)|\big(\displaystyle{\sum_{ d\le n\atop d|\ell}}\frac{1}{d^{2\s-1}}\big)  \Big),   & (\frac{1}{2}<\s<1),
\cr n^{1-2\s} \,\mathcal O_\s  \Big( \log n + \displaystyle{\sum_{\ell\in \Z_* }} |c_{f_{\s}}(\ell) |\,\# \{ d\le n\,: d|\ell \}   \Big), &  (0<\s\le  \frac{1}{2}).
\end{cases}
\end{equation*}
\end{theorem}

     In the series above,  the summand of order $\ell$ not only depends on estimates of M\"obius sums $\sum_{\l\le x} \frac{\m(\l)}{\l^{2\s}}$, for instance conditionally to RH($\a$), or using \eqref{mobiusweighted}, but  also on the divisors of $\ell$ which are  less than $n$, namely  on the arithmetical structure of the support  of the Fourier coefficient sequence. 
 \vskip 2pt In   Theorem  \ref{fp1},  the series
 $$\sum_{\ell\in \Z_* } c_{f_\s}(\ell)\Big(\sum_{d\le n\atop d|\ell}\frac{1}{d^{2\s-1}}\Big),$$
  cannot be estimated in general. The term in parenthesis 
  can be   close to $n^{2(1-\s)}$ (for those $\ell$  such that $d\le n\,\Rightarrow\, d|\ell$), or close to $2$ when the Fourier coefficients are supported by a sequence of numbers having few divisors, typically a sequence of primes. In this case we get for instance if $\s=\frac12$,
\begin{equation*}  \label{} 
F_{n,1/2}(f)
= \Big[\frac{n}{\zeta(2)} + \mathcal O(\log n)
\Big]\int_0^1f_{1/2}(x) \dd x
+     \mathcal O(1),
 \end{equation*}
 if   $\sum_{\ell\in \Z_* } |c_{f_{1/2}}(\ell) |$  converges. 
 
\begin{theorem}[$\s=1$]\label{fp4}    Let $f:[0,1]\to \R$ be such that ${f_1}\in \mathcal C$.
   {\rm (i)} Assume that the series $\sum_{\ell \in \Z} c_{f_1}(\ell)$ is convergent. Then,
\begin{align*}  
&F_{n,1}(f)
=\      \frac{\int_0^1f_1(x) \dd x }{\zeta(2)}\,\log n
     +
  \sum_{\ell\in \Z_* } c_{f_1}(\ell)  \sum_{ d\le n\atop d|\ell}\frac{1}{d}\Big(\sum_{\l\le \frac{n}{d}} \frac{\m(\l)}{\l^{2 }}\Big)
    +\mathcal O(1)
 .
\end{align*}
 \vskip 3 pt   {\rm (ii)}   If the series $\sum_{\ell\in \Z } |c_{f_1}(\ell)| $ is convergent, then  
\begin{align*} 
&F_{n,1}(f)
=\        \frac{\int_0^1f_1(x) \dd x }{\zeta(2)}\,\log n
 +
 \mathcal O\Big( \sum_{\ell\in \Z_* } |c_{f_1}(\ell)|\big(\sum_{ d\le n\atop d|\ell}\frac{1}{d}\big)\Big)   .
\end{align*}
In particular, if the series $\sum_{\ell\in \Z } |c_{f_1}(\ell)|\s_{-1}(\ell) $ is convergent, then 
\begin{align*} 
&F_{n,1}(f)
=\        \frac{\int_0^1f_1(x) \dd x }{\zeta(2)}\,\log n + \mathcal O ( 1)   .
\end{align*}
\end{theorem} 
\begin{remark}We recall   that by \cite{Gr}, $\s_{-1}(\ell)=\mathcal O( \log\log \ell)$,  thus the above estimate holds true under the mild condition 
$$\sum_{\ell\in \Z } |c_{f_1}(\ell)|\log\log \ell <\infty.$$ \end{remark}

The following Theorem concerns   Farey sums $F_n(f)$ (case $\s=0$) and provides a simple formula for the  error term $E_n(f)$ under   minimal assumption, as well as  a new estimate.  
 \begin{theorem}  \label{fc1}
   Let   ${f }\in \mathcal C$.
 {\rm (i)} If the series $\sum_{\ell\in \Z} c_f(\ell)$ is convergent, then
 \begin{eqnarray*}   E_n(f) 
 \, =\, \Big(f(1)-\sum_{\ell\in \Z} c_f(\ell)\Big)  + \sum_{\ell\in \Z_*} c_f(\ell) \Big(\sum_{ d\le n\atop d|\ell}d  M \big(\frac{n}{d} \big)\Big).
\end{eqnarray*}
 \vskip 3 pt \noi {\rm (ii)} If the series $\sum_{\ell\in \Z_*} |c_f(\ell)| d(\ell)$ is convergent, where $d(\ell)$ is the divisor function,  then  
 $$E_n(f)\, =\, o(n). $$
 In addition if  $f(1)   =\sum_{\ell\in \Z } c_{f}(\ell)$, then
\begin{eqnarray*}\sum_{n\ge 1}  \frac{|E_n(f)|
}{n^2} \, <\, \infty.
\end{eqnarray*}
 \end{theorem}

\begin{remarks}\rm  According to a classical estimate, $E_n(f)\, =\, o(n)$  if    $\sum_{\ell\in \Z_*} |c_f(\ell)|e^{c\log n / \log\log n}$ is convergent for some $c>\log 2$.  Codec\`a and Perelli \cite{CP} showed using Euler-McLaurin's formula and a lemma due to F\'ejer, that if $f$ is absolutely continuous, then $E_n(f)\, =\, o(n)$. See also Corollary p.\,105 in Mikolas \cite{M1}. This can in fact be improved if  in addition $f'\in L^p[0,1]$ for some $p\in (1,2]$, by using Vinogradov-Korobov's estimate. They further showed, under assumption RH($\a$), and using example  \ref{Jordan} below,  that for every $\e>0$, there exists an absolutely continuous function $f_\e:[0,1]\to \R$ such that 
$E_n(f)\, =\, \O_\e(n^{1-\e})$. 
\end{remarks}

 \begin{corollary}\label{corfcoeff} Assume that RH($\a$)  holds. Then for any  ${f }\in \mathcal C$ such that   the series 
 $$\sum_{\ell\in \Z} |c_f(\ell)| \s_{1-\a'}(\ell)$$
  converges for some $\a'>\a$, 
   we have
 $$E_n(f)  \, =\, \mathcal O(n^{\a'}).$$
 \end{corollary}
\begin{proof} As  RH($\a$) implies   $M(x)=\mathcal O_\e(x^{\a+\e})$, we have
\begin{eqnarray*}  
   \sum_{ d\le n\atop d|\ell}d\,\big|M \big(\frac{n}{d} \big)\big|\,\le\, C_{\a'}\,n^{\a'} \sum_{ d\le n\atop d|\ell}d^{1-\a'}\,\le\, C_{\a'}\,n^{\a'} \s_{1-\a'}(\ell).\end{eqnarray*}
 Thus
 \begin{eqnarray*}      \sum_{\ell\in \Z_*} |c_f(\ell)| \Big|\sum_{ d\le n\atop d|\ell}d  M \big(\frac{n}{d} \big)\Big| & \le &   C_{\a'}\,n^{\a'}  \sum_{\ell\in \Z_*} |c_f(\ell)| \s_{1-\a'}(\ell)\, \le \ C_{\a'}\,n^{\a'},
  \end{eqnarray*}     
 which by Theorem \ref{fc1} implies that $ E_n(f) 
 \, =\, \mathcal O_{\a'}(n^{\a'})$.
 \end{proof} 

\vskip 7pt 

 {\bf Discussion:} We 
 compare our assumptions  
  with the ones made  in \cite{CP}. First note, as a consequence of \eqref{farey.riemann} and of the fact that  $\sum_{ d\le n}  M \big(\frac{n}{d} \big) =1$  (see \eqref{3.12A}),    that
\begin{equation}\label{Enf1} E_n(f)
\ =\ \sum_{d\le n}  \tilde D_f(d) \,M\big(\frac{n}{d}\big),\qq \qq  \tilde D_f(n)=D_f(n)
- n\int_0^1 f(x)\dd x.
 \end{equation}
In \cite[Th.\,1\&2]{CP},
 assuming the validity of RH($\a$), it was proved that
$E_{n}(f)= \mathcal O_\e ( n^{(\a\vee {1}/{p})-\e})$ for all 
$f\in \mathcal C_{\rm cp}$, and the estimate is nearly optimal (\cite[Th.\,1\&2]{CP}).  
The authors underlined the   link between $\tilde D_f(n)$ and 
  Fourier coefficients  $(c_{f'}(m))_m$ of $f'$,   which follows from Euler-MacLaurin sum formula  (for absolutely continuous $f$, \cite[p.\,417]{CP}), 
  \begin{equation}\label{approx.term}\tilde D_f(n) \ =\  -\frac{1}{2\pi} \sum_{k=1}^\infty \frac{c_{f'}(kn)}{k},\qq\qq (f(0)= f(1)).\end{equation}
Note, however   (Bary \cite[Vol.\,I]{B}   p.\,78) that for  absolutely continuous $f$, $c_{f'_\s}(\ell)= i\ell c_{f_\s}(\ell)$, $\ell \in \Z_*$;     thus  the link with the   Fourier coefficients of $f$ is direct. 
\vskip 2pt
We  further  observe that  assumption  $f\in  \mathcal C_{\rm cp}$ implies that    
\begin{equation}\label{cp.ass}\sum_{\ell\in \Z_* }  |c_{f}(\ell)|\ell^\e<\infty,\qq\qq (\e< 1-1/p)\,
\end{equation} 
which is much stronger than our assumptions.
 Indeed, by  Hausdorff-Young's theorem,
\begin{eqnarray*}  
  \sum_{\ell\in \Z_* } |c_{f}(\ell)| \ell^\e& =& \sum_{\ell\in \Z_* } |c_{f'}(\ell)| \,\ell^{-(1-\e)}  
 \ \le \    \Big( \sum_{\ell\in \Z_* } |c_{f'}(\ell)|^q 
\Big)^{1/q}\Big( \sum_{\ell\in \Z_* } |\ell|^{- p(1-\e)}
\Big)^{1/p}
\cr  &\le & C \|{f'}\|_p\, \Big( \sum_{\ell\in \Z_* } |\ell|^{- p(1-\e)}
\Big)^{1/p}
\,<\,\infty
\, . \end{eqnarray*}

In the next Theorem we provide a sharp estimate of   the quadratic Farey sum $F_{n,\s}=F_{n,\s}(h)$, $h(x)=\cos 2\pi x$, recalling that $F_n(h)=M(n)$.
  \begin{theorem} \label{t1a}Let $0< \s<1$. Then
   \begin{eqnarray*}
 F_{n,\s}&=& \frac{ n^{2(1-\s)} } {2(1-\s)\zeta(2)}\int_{0}^{1} \frac{\cos 2\pi t}{t^\s}\dd t\, +  \mathcal O_\s\big(n^{1-\s}\big).
\end{eqnarray*}
  \end{theorem}
\vskip 1 pt   (The integral  term cannot be expressed elementarily.)

\vskip 5 pt 
 
Before passing to the proofs, let us give one more   example of Farey sums, linked to Euler's generalized totient function.
\begin{example}\label{Jordan}
Let    $\zeta_h (n) = n^h$  with summatory function $M_{\zeta_h}(x)$.         Let  $J_s(n)
 =\zeta_s*\m(n)=\sum_{d|n} d^s \m(n/d) $ be Euler's generalized totient function. 
  Let $f(x)= c(a)\sum_{\ell \in \Z^*} |\ell|^{-a} e_\ell(x)$,  $1<a<2$, where $c(a)=\big(\sum_{\ell \in \Z^*} |\ell|^{-a}\big)^{-1}$. 
  \vskip 2 pt Then 
  \begin{eqnarray} \label{jordan}F_n (f)&=&   \sum_{ k\le n} J_{1-a}(k),\end{eqnarray}
    the summatory function of   $J_{1-a}(n)$  therefore writes as a Farey sum. First note that $f(1)=f(0) =1
 $, and so by   Theorem \ref{fc1},
\begin{eqnarray*}   E_n(f)=F_n(f)& =&
  \sum_{\ell\in \Z_*} c_f(\ell) \Big(\sum_{ d\le n\atop d|\ell}d  M \big(\frac{n}{d} \big)\Big)\,=\,  c(a)\sum_{\ell\in \Z_*} |\ell|^{-a} \sum_{ d\le n\atop d|\ell}d \Big(\sum_{\l\le n/d}  \m(\l)\Big)
 \cr &=&c(a) \sum_{\l\le n}  \m(\l)\Big\{\sum_{ d\le n/\l}d\sum_{\ell\in \Z_*\atop d|\ell} |\ell|^{-a}   \Big\}\,=\,  \sum_{\l\le n}  \m(\l)\Big\{\sum_{ d\le n/\l}d^{1-a}   \Big\}.
\end{eqnarray*}

    By  well-known formula for partial sums of a Dirichlet product (\cite[Th.\,3.10]{A}),  
     
 \begin{eqnarray*} 
   \sum_{ k\le n}\m(k)\sum_{ d\le n/k} d^{1-a}
\,= \,\sum_{ k\le n}\m(k)M_{\zeta_{1-a}}(\frac{n}{k})
\, =\,  \sum_{ k\le n}\m* \zeta_{1-a}(k)= \sum_{ k\le n}J_{1-a}(k),
\end{eqnarray*}  
whence \eqref{jordan}.
\end{example}


\section{Proofs of Proposition \ref{Fp1} and Theorems \ref{fp1}, \ref{fp4}, \ref{fc1}.} \label{s1}
  We first establish some auxiliary lemmas and intermediate results.

 \subsection{Preliminary results.}
  
\begin{proposition}\label{25}Let  $f \in \mathcal C$.  Then,
\begin{eqnarray*}  
\ {\rm (i)}\qq\qq\qq\quad D_f(d)& =&  f(1)+  \sum_{\ell\in \Z} c_{f}(\ell)\e_\ell (d),\qq\qq\qq\qq\qq\qq\qq\qq\qq
\end{eqnarray*}
recalling  that $\e_\ell (d)= d-1$ if $d|\ell$, and $\e_\ell (d)= -1$ if $d\not|\ell$, by \eqref{eps}.
\vskip 4 pt 
\noi \ {\rm (ii)} If the series $ \sum_{\ell\in \Z } c_{f}(\ell)$ converges, then
\begin{eqnarray*}D_f(d) & =&  f(1)+  d\sum_{\ell\in \Z\atop d|\ell } c_{f}(\ell) -\sum_{\ell\in \Z} c_{f}(\ell).
\end{eqnarray*}\end{proposition}
 \begin{proof}Since $f\in \mathcal C$, by  Dini's test (Bary \cite[Vol.\,I]{B}   p.\,113, see also p.\,114) the Fourier series of ${f}$
$$\s_{f}(x)=\sum_{\nu\in \Z} c_{f}(\nu)e^{2{\rm i}\pi\nu x},$$
converges to ${f}(x)$ (i.e. the partials sums are converging to ${f}(x)$) for any $x\in \Q\cap ]0,1[$.
\vskip 2 pt
Thus 
\begin{eqnarray} \label{Festgen0}  
{f}\big(\frac{\k}{d}\big)\ =\   \sum_{\ell\in \Z} c_{f}(\ell)e^{2{\rm i}\pi\ell \frac{\k}{d}}, \qq \quad \k=1, 2, \ldots, d-1
\end{eqnarray}
whence 
\begin{eqnarray*}  
\sum_{\k=1}^{d-1} {f}\big(\frac{\k}{d}\big)\,=\, \sum_{\k=1}^{d-1}\sum_{\ell\in \Z} c_{f}(\ell)e^{2{\rm i}\pi\ell \frac{\k}{d}}
\,=\,    \sum_{\ell\in \Z} c_{f}(\ell)\big(\sum_{\k=1}^{d-1}e^{2{\rm i}\pi\ell \frac{\k}{d}}\big),
\end{eqnarray*}
  permutation between finitely many convergent series being permitted. 
  As $\sum_{\k=1}^{d-1}e^{2{\rm i}\pi\ell \frac{\k}{d}}=\e_\ell (d)$, we get 
\begin{eqnarray*}  
D_f(d)& =&  f(1)+    \sum_{\ell\in \Z} c_{f}(\ell)\e_\ell (d),
\end{eqnarray*}
Further  if the series $\sum_{\ell\in \Z} c_{f}(\ell)
$ converges, we can write
\begin{eqnarray*}  
D_f(d)& =&  f(1)+    \sum_{\ell\in \Z} c_{f}(\ell)(\e_\ell (d)+1)-\sum_{\ell\in \Z} c_{f}(\ell)
\cr & =&  f(1)+  d\sum_{\ell\in \Z\atop d|\ell } c_{f}(\ell) -\sum_{\ell\in \Z} c_{f}(\ell),
\end{eqnarray*}
which 
completes the proof.
\end{proof} 
\begin{corollary}\label{25a}Let  $f \in \mathcal C$. Assume that the series $ A=\sum_{\ell\in \Z  } |c_{f}(\ell)|d(\ell)$ converges. Then,
\begin{equation*}  \sum_{d=1}^\infty\frac{1}{d}\big|\tilde  D_f(d) -f(1)   +\sum_{\ell\in \Z } c_{f}(\ell)\big|
 \le   A .
\end{equation*}
\end{corollary}
\begin{proof} By   Proposition \ref{25}-(ii),  since $c_f(0)= \int_0^1 f(x)\dd x$,
\begin{eqnarray*}
D_f(d) & =&  f(1)+  d\int_0^1 f(x)\dd x+d\sum_{\ell\in \Z_*\atop d|\ell } c_{f}(\ell) -\sum_{\ell\in \Z} c_{f}(\ell),
\end{eqnarray*}
  By \eqref{Enf1},  
$$\tilde D_f(d)=  f(1)+d\sum_{\ell\in \Z_*\atop d|\ell } c_{f}(\ell) -\sum_{\ell\in \Z} c_{f}(\ell).$$
Thus
\begin{equation*}  
\sum_{d=1}^\infty\frac{1}{d}\big|\tilde  D_f(d) -f(1)   +\sum_{\ell\in \Z } c_{f}(\ell)\big| \ =\  \sum_{  d=1}^\infty     \Big| \sum_{\ell\in \Z_* \atop d|\ell  } c_{f}(\ell)\Big|   \ \le \    \sum_{\ell\in \Z_*  } |c_{f}(\ell)|\big( \sum_{  d|\ell} 1\big)= A,\end{equation*}
as claimed.\end{proof}
  
 \begin{lemma}\label{l2} Let $f:\Q\to \R$ be arbitrary and $\s\in \C$. Then, 
\begin{equation*} 
  F_{n,\s}(f)  \ =\ \sum_{d =1}^n\frac{1}{ d^{2\s}}\,\m*D_{f_\s}(d)\ =\ \sum_{d\le n} \frac{D_{f_\s}(d)}{d^{2\s}}
 \, \Big(\,\sum_{\l\le \frac{n}{d}} \, \frac{\m( \l)} { \l^{2\s}}\,\Big).
\end{equation*}
  \end{lemma}
\begin{proof}[Proof of Lemma \ref{l2}]  We recall that  
 \begin{equation}\label{inversemu}\sum_{d|n} \m(d)=\d(n)\qq\hbox{ where}\qq \d (n)=\begin{cases} 1  &\quad {\rm if}\ n=1, \cr
0&\quad {\rm   unless}.\end{cases}
\end{equation}
Then
\begin{eqnarray*}  
F_{n,\s}(f)&=& \sum_{\frac{k}{\ell} \in \F_{n}}\frac{1}{k^\s\ell^\s}f\big(\frac{k}{\ell}\big) 
\, =\,  \sum_{1\le k\le \ell \le n\atop (k,\ell) =1}\frac{1}{k^\s\ell^\s}f\big(\frac{k}{\ell}\big)=\sum_{\ell =1}^n\frac{1}{ \ell^\s}\sum_{k =1}^\ell\d\big((k,\ell)\big) \frac{1}{k^\s }f\big(\frac{k}{\ell}\big)
\cr &=& 
  \sum_{\ell =1}^n\frac{1}{ \ell^\s}\sum_{k =1}^\ell\Big(\sum_{u|(k,\ell)} \m(u)\Big)  \frac{1}{k^\s }f\big(\frac{k}{\ell}\big)
\, =\,  \sum_{\ell =1}^n\frac{1}{ \ell^\s}\sum_{u| \ell } \m(u) \sum_{k =1\atop u|k}^\ell   \frac{1}{k^\s }f\big(\frac{k}{\ell}\big).
\end{eqnarray*}
We write $k=\k u$ with $\k\le \frac{\ell}{u}$, and get 
\begin{eqnarray*}  
&= & \sum_{\ell =1}^n\frac{1}{ \ell^\s}\sum_{u| \ell } \frac{\m(u) }{u^\s} \sum_{1\le \k \le \frac{\ell}{u}}   \frac{1}{\k^\s }f\big( \frac{u\k}{\ell}\big).\qq \qq \qq \qq \qq \quad
\end{eqnarray*}
Now we write the divisors $u$ of $\ell$ under the form $u=\frac{\ell}{d}$,   $d$ running along all divisors of $\ell$, and   continue as follows
 \begin{eqnarray*} \qq  &= & \sum_{\ell =1}^n\frac{1}{ \ell^\s}\sum_{d| \ell } \frac{\m(\frac{\ell}{d}) }{(\frac{\ell}{d})^\s} \sum_{1\le \k \le d}   \frac{1}{\k^\s }f\big( \frac{ \k}{d}\big)
\, =\,\sum_{\ell =1}^n\frac{1}{ \ell^{2\s}}\sum_{d| \ell }  \m\big(\frac{\ell}{d}\big)   \sum_{1\le \k \le d}   \frac{1}{(\frac{ \k}{d})^\s }f\big( \frac{ \k}{d}\big)
\cr &=& \sum_{\ell =1}^n\frac{1}{ \ell^{2\s}}\sum_{d| \ell }  \m\big(\frac{\ell}{d}\big)   D_{f_\s}(d)
 \cr  &=& 
  \sum_{\ell =1}^n\frac{1}{ \ell^{2\s}}\,\m*D_{f_\s}(\ell)\, =\,  \sum_{d\le n} D_{f_\s}(d)
 \, \sum_{\ell\le n\atop d|\ell } \, \m(\frac{\ell}{d})  \frac{1}{ \ell^{2\s}}.
\end{eqnarray*}
Writing $\ell = \l d$ with $\l\le \frac{n}{d}$ in the last sum,  finally gives
 \begin{eqnarray*}  
  F_{n,\s}(f)  &=& \sum_{\ell =1}^n\frac{1}{ \ell^{2\s}}\,\m*D_{f_\s}(\ell)\,=\,\sum_{d\le n} \frac{D_{f_\s}(d)}{d^{2\s}}
 \, \Big(\,\sum_{\l\le \frac{n}{d}} \, \frac{\m( \l)} { \l^{2\s}}\,\Big).
\end{eqnarray*}
\end{proof}

We also need the following lemma.
\begin{lemma} \label{i-ii}   We have the following estimates.
\vskip 4 pt \noi   {\rm(a)} {\rm(}$\frac{1}{2}<\s<1${\rm)}
\begin{eqnarray*}
  {\rm (a\!\cdot\! 1)} & & 
  \sum_{ d\le n}\frac{1}{d^{2\s-1}} \big(\sum_{\l\le \frac{n}{d}} \frac{\m(\l)}{\l^{2\s}}\big)\ =\   \frac{n^{2(1-\s)}} {2(1-\s)\zeta(2)}+ \frac{\zeta(2\s-1)}{\zeta(2\s)}+ \mathcal O \big(n^{1-2\s}\log n \big),
\cr {\rm (a\!\cdot\! 2)} & & 
\ \  \, \sum_{ d\le n}\frac{1}{d^{2\s }} \big(\sum_{\l\le \frac{n}{d}} \frac{\m(\l)}{\l^{2\s}}\big)\ =\   1+ \mathcal O \big(n^{1-2\s}\log n \big),
   \end{eqnarray*}
recalling that 
$\zeta(s) = \lim_{x\to\infty} \big(\sum_{n\le x} \frac1{n^s}- \frac{x^{1-s}}{1-s}\big)$, $0<s<1$.
\vskip 2 pt \noi {\rm(b)}  {\rm(}$0<\s< \frac{1}{2}${\rm)} 
\begin{eqnarray*}
  {\rm (b\!\cdot\! 1)} & & 
  \sum_{ d\le n}\frac{1}{d^{2\s-1}} \big(\sum_{\l\le \frac{n}{d}} \frac{\m(\l)}{\l^{2\s}}\big)\ =\   \frac{n^{2(1-\s)}} {2(1-\s)\zeta(2)}+ \mathcal O \big(n^{1-2\s}\log n \big),
\cr {\rm (b\!\cdot\! 2)} & & 
\ \  \, \sum_{ d\le n}\frac{1}{d^{2\s }} \big(\sum_{\l\le \frac{n}{d}} \frac{\m(\l)}{\l^{2\s}}\big)\ =\   \mathcal O \big(n^{1-2\s}\big).
   \end{eqnarray*}

\vskip 2 pt \noi {\rm(c)}  If $\s= \frac12$, then
\begin{eqnarray*}
  {\rm (c\!\cdot\! 1)} & & 
  \sum_{ d\le n} \big(\sum_{\l\le \frac{n}{d}} \frac{\m(\l)}{\l}\big)\ = \  \frac{n}{\zeta(2)} + \mathcal O(\log n) ,
\cr {\rm (c\!\cdot\! 2)} & & 
\ \  \, \sum_{ d\le n}\frac{1}{d}\, \big|\sum_{\l\le \frac{n}{d}} \frac{\m(\l)}{\l}\big|\ \le \ B,\qq \quad n\ge 1,
\end{eqnarray*}where $B$ is an absolute constant.
 \vskip 2 pt \noi {\rm(d)}  If $\s=1$, then
\begin{eqnarray*}
  {\rm (d\!\cdot\! 1)} & & 
  \sum_{ d\le n} \frac{1}{d}\big(\sum_{\l\le \frac{n}{d}} \frac{\m(\l)}{\l^2}\big)\ = \  \frac{\log n}{\zeta(2)} + \mathcal O(1) ,
\cr {\rm (d\!\cdot\! 2)} & & 
\ \  \, \sum_{ d\le n}\frac{1}{d^2}\, \big|\sum_{\l\le \frac{n}{d}} \frac{\m(\l)}{\l^2}\big|\ = \  \mathcal O(1). 
\end{eqnarray*}
 \end{lemma} 
 \begin{proof}(i) 
We quote the estimates (\cite[Th.\,3.2]{A}),
\begin{equation}\label{re.zetasum}
\sum_{1\le \ell \le m} \ell^{1-2\s}=
\begin{cases}\frac{m(m+1)}{2} &\qq \hbox{if $\s=0$},
\cr
\frac{m^{2(1-\s)}}{2(1-\s)} + \mathcal O(m^{1-2\s}) &\qq \hbox{if $0<\s <1/2$},
\cr m &\qq \hbox{if $\s=1/2$},
\cr
\frac{m^{2(1-\s)}}{2(1-\s)} + \zeta(2\s-1) + \mathcal O\left(m^{1-2\s}\right) &\qq \hbox{if $\frac{1}{2}<\s<1$},
\cr
\log m + \g +\mathcal O(\frac1m) &\qq \hbox{if $\s=1$},
\end{cases}
\end{equation}
where $\g$ is Euler's constant 
 
\vskip 2 pt (a)  We get with \eqref{re.zetasum}, 
\begin{eqnarray*}  
 \sum_{ d\le n}\frac{1}{d^{2\s-1}} \Big[\sum_{\l\le \frac{n}{d}} \frac{\m(\l)}{\l^{2\s}}\Big]
 &=&  \sum_{\l\le n} \frac{\m(\l)}{\l^{2\s}}  \Big[\sum_{ d\le \frac{n}{\l}}\frac{1}{d^{2\s-1}}  \Big]
 \cr  &=&  \sum_{\l\le n} \frac{\m(\l)}{\l^{2\s}} \, \Big[\,\frac{1} {2(1-\s)} \big(\frac{n}{\l}\big)^{2(1-\s)} + \zeta(2\s-1) + \mathcal O \Big(\big(\frac{n}{\l} \big)^{1-2\s}  \Big) \Big]
\cr  &=& \frac{n^{2(1-\s)}} {2(1-\s)  }\Big(\sum_{\l\le n} 
\frac{\m(\l)}{\l^{2 }}\Big) + \zeta(2\s-1) \Big( \sum_{\l\le n} \frac{\m(\l)}{\l^{2\s}}  \Big)\cr & & \quad +  \mathcal O \Big(\,n^{1-2\s}\sum_{\l\le n} \frac{1}{\l}\,  \Big)
 \cr  &=&   \frac{n^{2(1-\s)}} {2(1-\s)  }\Big(\frac{1} { \zeta(2)}  + \mathcal O \big(n^{-1}  \big)\Big) + \zeta(2\s-1) \Big( \frac{1} { \zeta(2\s)}  + \mathcal O \big(n^{ 1-2\s}  \big) \Big)\cr & &\quad +  \mathcal O \big(n^{1-2\s}\log n \big)
 \cr &=& \frac{n^{2(1-\s)}} {2(1-\s)\zeta(2)}+ \frac{\zeta(2\s-1)}{\zeta(2\s)}+ \mathcal O \big(n^{1-2\s}\log n \big).  
  \end{eqnarray*}
Further,
 \begin{eqnarray*}  
 \sum_{ d\le n}\frac{1}{d^{2\s }} \Big[\sum_{\l\le \frac{n}{d}} \frac{\m(\l)}{\l^{2\s}}\Big]
 &=&  \sum_{\l\le n} \frac{\m(\l)}{\l^{2\s}}  \Big[\sum_{ d\le \frac{n}{\l}}\frac{1}{d^{2\s }}  \Big]
\cr &=&   \sum_{\l\le n} \frac{\m(\l)}{\l^{2\s}}  \Big[ \, \zeta(2\s) +\mathcal O \Big(\big(\frac{n}{\l} \big)^{1-2\s}    \Big)\,  \Big]
\cr &=&  \zeta(2\s) \sum_{\l\le n} \frac{\m(\l)}{\l^{2\s}}   +\mathcal O \Big(n^{1-2\s}\sum_{\l\le n} \frac{1}{\l}\,    \Big)
\cr &=&  \zeta(2\s) \Big(\frac{1}{\zeta(2\s)}+\mathcal O \big(n^{1-2\s}\big)\Big)    +\mathcal O \big(n^{1-2\s} \log n\,    \big)
\cr &=& 1   +\mathcal O \big(n^{1-2\s} \log n\,    \big).
  \end{eqnarray*}

  \vskip 2 pt (b) For ${\rm (b\!\cdot\! 1)}$ one gets similarly 
\begin{eqnarray*}  
 \sum_{ d\le n}\frac{1}{d^{2\s-1}} \Big[\sum_{\l\le \frac{n}{d}} \frac{\m(\l)}{\l^{2\s}}\Big]
 &=&  \sum_{\l\le n} \frac{\m(\l)}{\l^{2\s}}  \Big[\sum_{ d\le \frac{n}{\l}}\frac{1}{d^{2\s-1}}  \Big]
 \cr  &=&  \sum_{\l\le n} \frac{\m(\l)}{\l^{2\s}} \, \Big[\,\frac{1} {2(1-\s)} \big(\frac{n}{\l}\big)^{2(1-\s)} + \mathcal O \Big(\big(\frac{n}{\l} \big)^{1-2\s}  \Big) \Big]
\cr  &=& \frac{n^{2(1-\s)}} {2(1-\s)  }\Big(\sum_{\l\le n} 
\frac{\m(\l)}{\l^{2 }}\Big) +  \mathcal O \Big(\,n^{1-2\s}\sum_{\l\le n} \frac{1}{\l}\,  \Big)
\cr  &=& \frac{n^{2(1-\s)}} {2(1-\s)\zeta(2)}+  \mathcal O \big(n^{1-2\s}\log n \big) . \end{eqnarray*}
 Moreover,
 \begin{eqnarray*}  
 \sum_{ d\le n}\frac{1}{d^{2\s }} \Big[\sum_{\l\le \frac{n}{d}} \frac{\m(\l)}{\l^{2\s}}\Big]
 &=&  \sum_{\l\le n} \frac{\m(\l)}{\l^{2\s}}  \Big[\sum_{ d\le \frac{n}{\l}}\frac{1}{d^{2\s }}  \Big]
\cr &=&   \sum_{\l\le n} \frac{\m(\l)}{\l^{2\s}}  \Big[ \big(\frac{n}{\l} \big)^{ 1-2\s}\frac{1}{1-2\s} + \zeta(2\s) +\mathcal O \Big(\big(\frac{n}{\l} \big)^{ -2\s}    \Big)\,  \Big]
\cr &=&  \frac{n^{1-2\s}}{1-2\s}  \sum_{\l\le n} \frac{\m(\l)}{\l}    + \zeta(2\s)  \sum_{\l\le n} \frac{\m(\l)}{\l^{2\s}} +\mathcal O \big(n^{ -2\s} \sum_{\l\le n}1\big).
  \end{eqnarray*}
As \cite[Th.\,8.17]{BD}  $M(x)= \mathcal O( x\, e^{-c\sqrt{\log x}})$ for some positive number $c$,  it follows 
 by using Abel summation that
\begin{equation}\label{mobiusweighted}\sum_{ \l \le x}\frac{\m(\l)} {\l^{2\s}} =\mathcal O_\s (x^{1-2\s}\, e^{-c\sqrt{\log x}}),\qq \qq (0<\s<\frac12).
\end{equation}
 
  We deduce
 \begin{eqnarray*}  
 \sum_{ d\le n}\frac{1}{d^{2\s }} \Big[\sum_{\l\le \frac{n}{d}} \frac{\m(\l)}{\l^{2\s}}\Big]
 &=&   \mathcal O_\s (n^{1-2\s}).
  \end{eqnarray*}
   \vskip 2 pt (c) For ${\rm (c\!\cdot\! 1)}$, see Exercise 4-(b) in \cite[Ch.\,3]{A}. 
As to ${\rm (c\!\hskip -1 pt\cdot\!\hskip -1 pt2)}$, 
  this  is Rubel's estimate (\cite[Prob.\,8.20]{BD}).
 \vskip 2 pt (d) We have 
 \begin{eqnarray*}  
 \sum_{ d\le n}\frac{1}{d} 
 \Big[\sum_{\l\le \frac{n}{d}} \frac{\m(\l)}{\l^{2 }}\Big]
 &=& \sum_{ d\le n}\frac{1}{d}\Big[\frac{1}{\zeta(2 )} +\mathcal O( 
 \frac{d}{n})\Big]
\,=\,  \frac{\log n}{\zeta(2 )} +\mathcal O( 
 1),
  \end{eqnarray*}
 and obviously, 
 \begin{eqnarray*}  
 \sum_{ d\le n}\frac{1}{d^{2  }} 
 \Big[\sum_{\l\le \frac{n}{d}} \frac{\m(\l)}{\l^{2 }}\Big]
 &=&  \mathcal O( 
 1).
  \end{eqnarray*}
 \end{proof}

 \begin{proof}[Proof of Proposition \ref{Fp1}]
 (i)   Using Lemma \ref{l2} we have
\begin{eqnarray*}  
& &F_{n,\s}(f)\ =\ \sum_{ d\le n}\frac{D_{f_\s}(d)}{d^{2\s}}\Big(\sum_{\l\le \frac{n}{d}} \frac{\m(\l)}{\l^{2\s}}\Big)
\cr &=& {f_\s}(1)\Big(\sum_{ d\le n}\frac{1}{d^{2\s}} \big(\sum_{\l\le \frac{n}{d}} \frac{\m(\l)}{\l^{2\s}}\big)\Big)
+\sum_{ d\le n}\frac{1}{d^{2\s}}\Big(\sum_{\k=1}^{d-1}{f_\s}\big(\frac{\k}{d}\big)\Big)\Big(\sum_{\l\le \frac{n}{d}} \frac{\m(\l)}{\l^{2\s}}\Big),
\end{eqnarray*}
where we isolated   the  term related to $ {f_\s}(1)$.
 As $f_\s \in \mathcal C$, by Proposition \ref{25},
 \begin{eqnarray*}  
\sum_{\k=1}^{d-1}{f_\s}\big(\frac{\k}{d}\big) & =& D_{f_\s}(d)- f_\s(1)\,=\,  (d-1) c_{f}(0)+  \sum_{\ell\in \Z_*} c_{f}(\ell)\e_\ell (d).
\end{eqnarray*}
Thus,
\begin{eqnarray*}
F_{n,\s}(f)
 &=& {f_\s}(1)\Big(\sum_{ d\le n}\frac{1}{d^{2\s}} \big(\sum_{\l\le \frac{n}{d}} \frac{\m(\l)}{\l^{2\s}}\big)\Big)
+ c_{f}(0)\sum_{ d\le n}\frac{d-1}{d^{2\s}}\Big(\sum_{\l\le \frac{n}{d}} \frac{\m(\l)}{\l^{2\s}}\Big)
\cr && +\sum_{ d\le n}\frac{1}{d^{2\s}}\Big(\sum_{\ell\in \Z_*} c_{f}(\ell)\e_\ell (d)\Big)\Big(\sum_{\l\le \frac{n}{d}} \frac{\m(\l)}{\l^{2\s}}\Big)
\cr &=  & {f_\s}(1)\Big(\sum_{ d\le n}\frac{1}{d^{2\s}} \big(\sum_{\l\le \frac{n}{d}} \frac{\m(\l)}{\l^{2\s}}\big)\Big)
+c_{f_\s}(0)\Big[\sum_{ d\le n}\frac{d-1}{d^{2\s}} \Big(\sum_{\l\le \frac{n}{d}} \frac{\m(\l)}{\l^{2\s}}\Big)
\Big]
   \cr & &\qq+
   \sum_{\ell\in \Z_*} c_{f_\s}(\ell) \Big[\sum_{ d\le n}\frac{1}{d^{2\s}}    \e_\ell(d)  \Big(\sum_{\l\le \frac{n}{d}} \frac{\m(\l)}{\l^{2\s}}\Big)\Big]\end{eqnarray*}
since 
only a finite number of convergent series is involved.
\vskip 3 pt
(ii) It is easy to observe  with Proposition \ref{Fp1}
 (applied with  $\s=0$ and thus ${f_\s}=f$)  that      we also have (recalling that $\Phi(n)=\sum_{ d\le n}d M \big(\frac{n}{d} \big)$),
\begin{eqnarray*}  
F_n(f) 
&=& f(1)\Big(\sum_{ d\le n} M \big(\frac{n}{d} \big)\Big)
+\Phi(n)
 \int_0^1f(x)\dd x 
  +
   \sum_{\ell\in \Z_*} c_f(\ell) \Big[\sum_{ d\le n}\e_\ell(d)  M \big(\frac{n}{d} \big)\Big],\end{eqnarray*}
where  $\e_\ell (d)$ is defined in \eqref{eps}.
 By Theorem 3.12 in \cite{A}, we have 
 \begin{eqnarray}\label{3.12A}
 \sum_{ d\le n}  M \big(\frac{n}{d} \big)=\sum_{u \le n} \m(u )\Big\lfloor \frac{n}{u }\Big\rfloor  =1.
 \end{eqnarray} So that
\begin{eqnarray*}  
F_n(f) 
&=& f(1)+\Phi(n)
 \int_0^1f(x)\dd x+   \sum_{\ell\in \Z_*} c_f(\ell) \Big[\sum_{ d\le n}\e_\ell(d)  M \big(\frac{n}{d} \big)\Big].\end{eqnarray*} \vskip 2 pt
Applying this to $\bar f= f-\int_0^1f(x)\dd x$ in place of $f$ gives ($c_{\bar f}(\ell)=c_{ f}(\ell )$, $\ell\neq 0$ and $c_{\bar f}(0)=0$).
\begin{eqnarray*} 
E_n(f)&=&F_n(f) -\Phi(n)\Big(\int_0^1f(x)\dd x\Big)
\cr &=& \Big(f(1)-\int_0^1f(x)\dd x\Big)+
   \sum_{\ell\in \Z_*} c_f(\ell)  \Big[\sum_{ d\le n}\e_\ell(d)  M \big(
   \frac{n}{d} \big)\Big] .
\end{eqnarray*}
\end{proof}
\subsection{Proof of Theorem \ref{fp1}}
   By Proposition \ref{Fp1},
  \begin{align*} 
F_{n,\s}(f)
=\ & {f_\s}(1)\Big(\sum_{ d\le n}\frac{1}{d^{2\s}} \big(\sum_{\l\le \frac{n}{d}} \frac{\m(\l)}{\l^{2\s}}\big)\Big)
+c_{f_\s}(0)\Big[\sum_{ d\le n}\frac{d-1}{d^{2\s}} \Big(\sum_{\l\le \frac{n}{d}} \frac{\m(\l)}{\l^{2\s}}\Big)
\Big]
   \cr &\qq+
   \sum_{\ell\in \Z_*} c_{f_\s}(\ell) \Big[\sum_{ d\le n}\frac{1}{d^{2\s}}    \e_\ell(d)  \Big(\sum_{\l\le \frac{n}{d}} \frac{\m(\l)}{\l^{2\s}}\Big)\Big].
\end{align*}
  Since  the series $\sum_{\ell \in \Z}c_{f_\s}(\ell)$ is convergent, we note that
  \begin{align*} 
  \sum_{\ell\in \Z_*} c_{f_\s}(\ell) &\Big[\sum_{ d\le n}\frac{1}{d^{2\s}}    \e_\ell(d)  \Big(\sum_{\l\le \frac{n}{d}} \frac{\m(\l)}{\l^{2\s}}\Big)\Big]
\cr &=      \sum_{\ell\in \Z_* } c_{f_\s}(\ell)  \sum_{ d\le n\atop d|\ell}\frac{1}{d^{2\s-1}}\Big(\sum_{\l\le \frac{n}{d}} \frac{\m(\l)}{\l^{2\s}}\Big)
-\Big(\sum_{\ell\in \Z_* } c_{f_\s}(\ell) \Big)   .
\end{align*}
So that
 \begin{align}  \label{conv.eps}
F_{n,\s}(f)
=\ & {f_\s}(1)\Big(\sum_{ d\le n}\frac{1}{d^{2\s}} \big(\sum_{\l\le \frac{n}{d}} \frac{\m(\l)}{\l^{2\s}}\big)\Big)
+c_{f_\s}(0)\Big[\sum_{ d\le n}\frac{d-1}{d^{2\s}} \Big(\sum_{\l\le \frac{n}{d}} \frac{\m(\l)}{\l^{2\s}}\Big)
\Big]
   \cr &\qq+
  \sum_{\ell\in \Z_* } c_{f_\s}(\ell)  \sum_{ d\le n\atop d|\ell}\frac{1}{d^{2\s-1}}\Big(\sum_{\l\le \frac{n}{d}} \frac{\m(\l)}{\l^{2\s}}\Big)
-\Big(\sum_{\ell\in \Z_* } c_{f_\s}(\ell) \Big)  .
\end{align}
\vskip 3 pt ({\bf a})  ($\frac{1}{2}<\s<1$) 
\vskip 2 pt \noi  (i) By applying 
Lemma \ref{i-ii}, we obtain
\begin{eqnarray*}  
& & F_{n,\s}(f) \ = \   {f_\s}(1)+c_{f_\s}(0)  \Big(  \frac{n^{2(1-\s)}} {2(1-\s)\zeta(2)}  + \frac{\zeta(2\s-1)}{\zeta(2\s)} -1\Big) 
\cr & & +   \sum_{\ell\in \Z_* } c_{f_\s}(\ell)  \sum_{ d\le n\atop d|\ell}\frac{1}{d^{2\s-1}}\Big(\sum_{\l\le \frac{n}{d}} \frac{\m(\l)}{\l^{2\s}}\Big)
-\Big(\sum_{\ell\in \Z_* } c_{f_\s}(\ell) \Big)  + \mathcal O (n^{1-2\s}\log n   )
\cr   &=&      \Big( \frac{\int_0^1f_\s(x) \dd x} {2(1-\s)\zeta(2)}  \Big)   n^{2(1-\s)}+   \sum_{\ell\in \Z_* } c_{f_\s}(\ell)  \sum_{ d\le n\atop d|\ell}\frac{1}{d^{2\s-1}}\Big(\sum_{\l\le \frac{n}{d}} \frac{\m(\l)}{\l^{2\s}}\Big)
    +   A    
    \cr & &+ \mathcal O (n^{1-2\s}\log n   )  ,\end{eqnarray*}
recalling that $A={f_\s}(1)+  \big(\int_0^1f_\s(x)\dd x\big)\big(  \frac{\zeta(2\s-1)}{\zeta(2\s)} -1\big)  -\sum_{\ell\in \Z_* } c_{f_\s}(\ell)$.
\vskip 3 pt  (ii) If  the series $\sum_{\ell \in \Z}|c_{f_\s}(\ell)|$ is convergent, then
\begin{eqnarray*}   \sum_{ d\le n\atop d|\ell}\frac{1}{d^{2\s-1}} 
\Big(\sum_{\l\le \frac{n}{d}} \frac{\m(\l)}{\l^{2\s}}\Big)
&=& \sum_{ d\le n\atop d|\ell}\frac{1}{d^{2\s-1}} \Big[\frac{1}{\zeta(2\s)}+\mathcal O\Big(\big( \frac{n}{d}\big)^{1-2\s}\Big) \Big]
\cr &=& \frac{1}{\zeta(2\s)} \sum_{ d\le n\atop d|\ell}\frac{1}{d^{2\s-1}} +\mathcal O\Big(n^{1-2\s}\#\{ d\le n: d|\ell\}\Big). 
\end{eqnarray*}  
Thus
\begin{eqnarray*}  
F_{n,\s}(f) &=&  {f_\s}(1)+c_{f_\s}(0)  \Big(  \frac{n^{2(1-\s)}} {2(1-\s)\zeta(2)}  + \frac{\zeta(2\s-1)}{\zeta(2\s)} -1\Big) -\Big(\sum_{\ell\in \Z_* } c_{f_\s}(\ell) \Big) 
\cr & & +   \frac{1}{\zeta(2\s)} \sum_{\ell\in \Z_* } c_{f_\s}(\ell)  \sum_{ d\le n\atop d|\ell}\frac{1}{d^{2\s-1}}  + \mathcal O (n^{1-2\s}\log n   ) 
\cr & &+ \mathcal O\Big( n^{1-2\s}\sum_{\ell\in \Z_* } |c_{f_\s}(\ell)|\#\{ d\le n: d|\ell\}\Big)
\cr &=&    \Big( \frac{\int_0^1f_\s(x) \dd x} {2(1-\s)\zeta(2)}  \Big) \, n^{2(1-\s)} +   \frac{1}{\zeta(2\s)} \sum_{\ell\in \Z_* } c_{f_\s}(\ell)\Big(\sum_{d\le n\atop d|\ell}\frac{1}{d^{2\s-1}}\Big) +A  \cr 
  & &   + n^{1-2\s}\, \mathcal O_\s  \Big( \log n+   \sum_{\ell\in \Z_* } |c_{f_{\s}}(\ell) |\,\# \{ d\le n\,: d|\ell \}   \Big)  .
\end{eqnarray*}

\vskip 3 pt ({\bf b})  ($0<\s< \frac{1}{2}$) 
\vskip 2 pt \noi  (i)
From \eqref{conv.eps} and Lemma \ref{i-ii} follows that  
 \begin{eqnarray*}
F_{n,\s}(f)
 &= & {f_\s}(1)\Big(\sum_{ d\le n}\frac{1}{d^{2\s}} \big(\sum_{\l\le \frac{n}{d}} \frac{\m(\l)}{\l^{2\s}}\big)\Big)
+c_{f_\s}(0)\Big[\sum_{ d\le n}\frac{d-1}{d^{2\s}} \Big(\sum_{\l\le \frac{n}{d}} \frac{\m(\l)}{\l^{2\s}}\Big)
\Big]
   \cr & &\qq+
  \sum_{\ell\in \Z_* } c_{f_\s}(\ell)  \sum_{ d\le n\atop d|\ell}\frac{1}{d^{2\s-1}}\Big(\sum_{\l\le \frac{n}{d}} \frac{\m(\l)}{\l^{2\s}}\Big)
-\Big(\sum_{\ell\in \Z_* } c_{f_\s}(\ell) \Big) 
\cr &=& \big({f_\s}(1)-c_{f_\s}(0)\big)\Big(\sum_{ d\le n}\frac{1}{d^{2\s}} \big(\sum_{\l\le \frac{n}{d}} \frac{\m(\l)}{\l^{2\s}}\big)\Big)
+c_{f_\s}(0)\Big[\sum_{ d\le n}\frac{1}{d^{2\s-1}} \Big(\sum_{\l\le \frac{n}{d}} \frac{\m(\l)}{\l^{2\s}}\Big)
\Big]
   \cr & & \qq+
  \sum_{\ell\in \Z_* } c_{f_\s}(\ell)  \sum_{ d\le n\atop d|\ell}\frac{1}{d^{2\s-1}}\Big(\sum_{\l\le \frac{n}{d}} \frac{\m(\l)}{\l^{2\s}}\Big)
-\Big(\sum_{\ell\in \Z_* } c_{f_\s}(\ell) \Big)
\cr &=& \mathcal O\big( n^{1-2\s}\big)+c_{f_\s}(0)\Big[ \frac{n^{2(1-\s)}}{2(1-\s)\zeta(2)}+ \mathcal O\big( n^{1-2\s}\log n\big)
\Big]
   \cr & &\qq+
  \sum_{\ell\in \Z_* } c_{f_\s}(\ell)  \sum_{ d\le n\atop d|\ell}\frac{1}{d^{2\s-1}}\Big(\sum_{\l\le \frac{n}{d}} \frac{\m(\l)}{\l^{2\s}}\Big)
-\Big(\sum_{\ell\in \Z_* } c_{f_\s}(\ell) \Big)
\cr &=& c_{f_\s}(0)\Big[ \frac{n^{2(1-\s)}}{2(1-\s)\zeta(2)}
\Big]+
  \sum_{\ell\in \Z_* } c_{f_\s}(\ell)  \sum_{ d\le n\atop d|\ell}\frac{1}{d^{2\s-1}}\Big(\sum_{\l\le \frac{n}{d}} \frac{\m(\l)}{\l^{2\s}}\Big)
   \cr & &\qq+ \mathcal O\big( n^{1-2\s}\log n\big)-\Big(\sum_{\ell\in \Z_* } c_{f_\s}(\ell) \Big) .
\end{eqnarray*}
\vskip 2 pt \noi (ii) If  the series $\sum_{\ell \in \Z}|c_{f_\s}(\ell)|$ is convergent, then
by \eqref{mobiusweighted},
\begin{eqnarray*}\sum_{ d\le n\atop d|\ell}\frac{1}{d^{2\s-1}}\Big|\sum_{\l\le \frac{n}{d}} \frac{\m(\l)}{\l^{2\s}}\Big|& =&\mathcal O_\s \Big(\sum_{ d\le n\atop d|\ell}\frac{1}{d^{2\s-1}}\big(\frac{n}{d}\big)^{1-2\s}\, e^{-c\sqrt{\log (\frac{n}{d})}}\Big)\cr &=&\mathcal O_\s \big(n^{1-2\s}\#\{d\le n: d|\ell\}\big).
\end{eqnarray*}
And so,
 \begin{align*}
F_{n,\s}(f)
=\ & c_{f_\s}(0)\Big[ \frac{n^{2(1-\s)}}{2(1-\s)\zeta(2)}
\Big]+n^{1-2\s} \, \mathcal O\Big( \log n+\sum_{\ell\in \Z_* } |c_{f_\s}(\ell) |\#\{d\le n: d|\ell\}\Big) .
\end{align*}

\vskip 3 pt ({\bf c})  ($\s= \frac{1}{2}$) 
\vskip 2 pt \noi  (i)
 Assume that the series $\sum_{\ell \in \Z}c_{f_{1/2}}(\ell)$ is convergent. By \eqref{conv.eps}, 

 \begin{eqnarray}  \label{} 
F_{n,1/2}(f)
&=&  {f_{1/2}}(1)\Big(\sum_{ d\le n}\frac{1}{d} \big(\sum_{\l\le \frac{n}{d}} \frac{\m(\l)}{\l}\big)\Big) 
 \cr & & +c_{f_{1/2}}(0)\Big[\sum_{ d\le n}\frac{d-1}{d} \Big(\sum_{\l\le \frac{n}{d}} \frac{\m(\l)}{\l}\Big)
\Big]
  \cr & &   \qq  +     \sum_{\ell\in \Z_* } c_{f_{1/2}}(\ell)  \sum_{ d\le n\atop d|\ell}\Big(\sum_{\l\le \frac{n}{d}} \frac{\m(\l)}{\l}\Big)
-\Big(\sum_{\ell\in \Z_* } c_{f_{1/2}}(\ell) \Big) .
\end{eqnarray}
 By Lemma \ref{i-ii} we first get
\begin{eqnarray*}
\sum_{ d\le n}\frac{d-1}{d} \Big(\sum_{\l\le \frac{n}{d}} \frac{\m(\l)}{\l}\Big)
&=& \sum_{ d\le n} \Big(\sum_{\l\le \frac{n}{d}} \frac{\m(\l)}{\l}\Big)
-\sum_{ d\le n}\frac{1}{d} \Big(\sum_{\l\le \frac{n}{d}} \frac{\m(\l)}{\l}\Big)
\cr &=& \frac{n}{\zeta(2)} + \mathcal O(\log n),\end{eqnarray*}
  next    
\begin{eqnarray*}  
F_{n,1/2}(f)
&=&  c_{f_{1/2}}(0)\Big[\frac{n}{\zeta(2)} + \mathcal O(\log n)
\Big]+     \sum_{\ell\in \Z_* } c_{f_{1/2}}(\ell)  \sum_{ d\le n\atop d|\ell}\Big(\sum_{\l\le \frac{n}{d}} \frac{\m(\l)}{\l}\Big)
  \cr & &   \qq  -\Big(\sum_{\ell\in \Z_* } c_{f_{1/2}}(\ell) \Big).
\end{eqnarray*}

\vskip 2 pt \noi (ii) By Theorem 3.13 in \cite{A},
\begin{eqnarray}  \label{A3.13} \Big| \sum_{\l\le x} \frac{\m(\l)}{\l}\Big| \le 1,
\end{eqnarray}
 if $x\ge 1$.
 Thus 
 \begin{eqnarray*}  
     \sum_{\ell\in \Z_* } |c_{f_{1/2}}(\ell) | \sum_{ d\le n\atop d|\ell}\Big| \sum_{\l\le \frac{n}{d}} \frac{\m(\l)}{\l}\Big|&\le &  \sum_{\ell\in \Z_* } |c_{f_{1/2}}(\ell) | \,\#\big\{ d\le n\,: d|\ell\big\},
       \end{eqnarray*}
 the series $\sum_{\ell\in \Z_* } |c_{f_{1/2}}(\ell) |$ converging by assumption. We therefore get
\begin{equation*}  \label{} 
\sum_{\frac{k}{\ell} \in \F_n}\frac{1}{(k\ell)^{1/2}} f\big(\frac{k}{\ell}\big)
= c_{f_{1/2}}(0)\Big[\frac{n}{\zeta(2)} + \mathcal O(\log n)
\Big]
+     \mathcal O \Big(\sum_{\ell\in \Z_* } |c_{f_{1/2}}(\ell) | \,\#\big\{ d\le n\,: d|\ell\big\}\Big).
 \end{equation*}
As $c_{f_{1/2}}(0)=\int_0^1f_{1/2}(x) \dd x$, this achieves the proof.


\subsection{Proof of Theorem \ref{fp4}}
If the series $\sum_{\ell\in \Z} c_f(\ell)$ converges, it follows from \eqref{conv.eps}   that
\begin{align}  \label{conv.eps1}
&F_{n,1}(f)
=\  {f_1}(1)\Big(\sum_{ d\le n}\frac{1}{d^{2 }} \big(\sum_{\l\le \frac{n}{d}} \frac{\m(\l)}{\l^{2}}\big)\Big)
+c_{f_1}(0)\Big[\sum_{ d\le n}\frac{1}{d} \Big(\sum_{\l\le \frac{n}{d}} \frac{\m(\l)}{\l^{2}}\Big)
\Big]
   \cr &-c_{f_1}(0)\Big[\sum_{ d\le n}\frac{ 1}{d^{2}} \Big(\sum_{\l\le \frac{n}{d}} \frac{\m(\l)}{\l^{2}}\Big)
\Big]+
  \sum_{\ell\in \Z_* } c_{f_1}(\ell)  \sum_{ d\le n\atop d|\ell}\frac{1}{d}\Big(\sum_{\l\le \frac{n}{d}} \frac{\m(\l)}{\l^{2 }}\Big)
\    -\Big(\sum_{\ell\in \Z_* } c_{f_1}(\ell) \Big)
    \cr  =\ &  c_{f_1}(0) \frac{\log n}{\zeta(2)}
     +
  \sum_{\ell\in \Z_* } c_{f_1}(\ell)  \sum_{ d\le n\atop d|\ell}\frac{1}{d}\Big(\sum_{\l\le \frac{n}{d}} \frac{\m(\l)}{\l^{2 }}\Big)
    +\mathcal O(1)
 .
\end{align}

And if   the series $\sum_{\ell\in \Z} |c_f(\ell)|$ converges, then
\begin{align}  \label{conv.eps2}
&F_{n,1}(f)
=\      c_{f_1}(0) \frac{\log n}{\zeta(2)}
     +
 \mathcal O\Big( \sum_{\ell\in \Z_* } |c_{f_1}(\ell)|\big(\sum_{ d\le n\atop d|\ell}\frac{1}{d}\big)\Big)   .
\end{align}


 \subsection{Proof of Theorem \ref{fc1}}
     (i) Since $$\sum_{ d\le n}\e_\ell(d)  M \big(\frac{n}{d} \big)=
\sum_{ d\le n\atop 
d|\ell}d  M \big(\frac{n}{d} \big)-\sum_{ d\le n}  M \big(\frac{n}{d} \big),$$
and by assumption the series $\sum_{\ell\in \Z} c_f(\ell)$ is convergent,  we get using Proposition \ref{Fp1}-(ii) and \eqref{3.12A},
\begin{eqnarray*}   E_n(f) &=&  f(1)-\int_0^1f(x)\dd x +
   \sum_{\ell\in \Z_*} c_f(\ell) \Big(\sum_{ d\le n\atop 
d|\ell}d  M \big(\frac{n}{d} \big)\Big)
 - \Big( \sum_{\ell\in \Z_*} c_f(\ell)\Big)\Big(\sum_{ d\le n}  M \big(\frac{n}{d} \big)\Big)
\cr &=&  f(1)-\sum_{\ell\in \Z} c_f(\ell) +
   \sum_{\ell\in \Z_*} c_f(\ell) \Big(\sum_{ d\le n\atop 
d|\ell}d  M \big(\frac{n}{d} \big)\Big)
.
\end{eqnarray*}

(ii) As $M(x)=o(x)$, we deduce from (i) and the assumption made that
 \begin{eqnarray*}   E_n(f)
 &=&  f(1)-\sum_{\ell\in \Z} c_f(\ell) +
 o\Big(n \sum_{\ell\in \Z_*} |c_f(\ell)| d(\ell) \Big)\ =\ o(n).
\end{eqnarray*}
 
Now if further  $f(1)   =\sum_{\ell\in \Z } c_{f}(\ell)$,   Corollary \ref{25a} reads 
\begin{equation}  \label{25ab}
 \sum_{d=1}^\infty\frac{1}{d}\big|\tilde  D_f(d) \big|
 \le   A ,
\end{equation}
recalling that $\tilde D_f(d)=D_f(d)
- d\int_0^1 f(x)\dd x$. By Theorem 8.17 of \cite{BD}, there exists a positive number $c$ such that
$$M(x) =\mathcal O(x e^{-c\sqrt{\log x}}).$$
By \eqref{Enf1}, $E_n(f)
\ =\ \sum_{d\le n}  \tilde D_f(d) \,M\big(\frac{n}{d}\big)$. Let $\varpi (n)$ be non-decreasing. Then,

\begin{eqnarray*}\sum_{n\ge 1}  \frac{|E_n(f)|
}{\varpi (n)}
&\le & \sum_{n\ge 1}  \frac{1
}{\varpi (n)}\sum_{d\le n}  |\tilde D_f(d) |\,\big| M\big(\frac{n}{d}\big)\big|
\cr 
&\le & C \sum_{n\ge 1}  \frac{n
}{\varpi (n)}\sum_{d\le n}  \frac{|\tilde D_f(d) |}{d} \,e^{-c\sqrt{\log \frac{n}{d}}} 
\cr 
&= &C \sum_{d\ge 1 } \frac{|\tilde D_f(d) |}{d}\sum_{n\ge d}  \frac{n
}{\varpi (n)}    \,e^{-c\sqrt{\log \frac{n}{d}}} 
\cr &\le  &C \sum_{d\ge 1 } \frac{|\tilde D_f(d) |}{d}\int_{  d}^\infty  \frac{t
}{\varpi (t)}    \,e^{-c\sqrt{\log \frac{t}{d}}} \dd t 
\cr &\le  &C \sum_{d\ge 1 }  |\tilde D_f(d) | \int_{  1}^\infty  \frac{u
}{\varpi (ud)}    \,e^{-c\sqrt{\log  u}} \dd u 
.
\end{eqnarray*}
 Taking $\varpi (v)= v^2$ we get, 
\begin{eqnarray*}\sum_{n\ge 1}  \frac{|E_n(f)|
}{n^2}
&\le &C \sum_{d\ge 1 } \frac{ |\tilde D_f(d) |}{d} \int_{  1}^\infty  \frac{1
}{u}    \,e^{-c\sqrt{\log  u}} \dd u \, <\, \infty,
\end{eqnarray*}
by   \eqref{25ab}.

 \section{Proof of Theorems \ref{fnfomega} and  \ref{t1a}.}
    \label{proof{t1a}}
\begin{proof}[{\bf Proof of Theorem \ref{fnfomega}.}] Let $0<Y_k\uparrow \infty$ with $k$, be an   integer sequence such that 
$$ {Y_{k+1}}/{e^{5 (\log\log Y_{k+1})^{3/2}}}\, \ge \,Y_k+2 , \qq \qq k= 1, 2, \ldots$$

We use the following precise result of Pintz \cite{Pi}: for $Y>c$ effective, there exists $x',x''\in [Ye^{-5 (\log\log Y)^{3/2}},Y]$ such that 
\begin{equation} \label{Pintz}
M(x')<-\frac{\sqrt{ x'}}{L}, \qq M(x'')>\frac{\sqrt{ x''}}{L}.
\end{equation}
$L=136000$. Assume that $Y_1>c$. Let $n= Y_{k+1}$. There thus  exist integers $m_j$,  $1\le j\le k$, such that 
$$m_j\in [Y_{j+1}\,e^{-5 (\log\log Y_{j+1})^{3/2}},Y_{j+1}], \qq {\rm and}\quad M(m_j)>\frac{\sqrt{ m_j}}{L}.$$
 For each $1\le j\le k$ let $d_j$ integer be such that $m_j=\lfloor\frac{n}{d_j}\rfloor$. 
  The numbers $d_j$ are mutually distinct since 
\begin{eqnarray*}\frac{d_{j-1}-d_j}{Y_{k+1}}&=&\frac{1}{\frac{Y_{k+1}}{d_{j-1}}}-\frac{1}{\frac{Y_{k+1}}{d_{j}}}
\,\ge\,\frac{1}{\lfloor\frac{Y_{k+1}}{d_{j-1}}\rfloor+1}-\frac{1}{\lfloor\frac{Y_{k+1}}{d_{j}}\rfloor} \,=\,  \frac{1}{m_{j-1}+1}-\frac{1}{m_j} 
\cr &\ge& \frac{1}{Y_j+1}-\frac{1}{Y_{j}+2} \,=\,\frac{1}{(Y_j+1)(Y_j+2)}\,>\, 0.
\end{eqnarray*}

Next we use the infinite M\"obius inversion formula which we recall.\vskip 3 pt
\noi{\it Infinite M\"obius inversion formula.} By Theorem 270 of \cite{HW}, we have the other (infinite) M\"obius inversion formula
\begin{equation}\label{HW} 
g(x) =\sum_{m=1}^\infty f(mx)\qq \Leftrightarrow \qq f(x) =\sum_{n=1}^\infty \m(n)g(nx),
\end{equation}
if for instance 
\begin{equation*} 
 \sum_{m,n=1}^\infty |f(mnx)| =\sum_{k=1}^\infty d(k) |f(kx)|<\infty.
\end{equation*}
It also suffices that
\begin{equation} \label{HWcondition}
\sum_{\nu=1}^\infty d(\nu)  |g(\nu x)|<\infty.
\end{equation}
As no proof is given in \cite{HW}, we provide it  for  sake of completeness.  First, under this assumption $f$ is obviously well defined. Now 
\begin{eqnarray*}\sum_{m=1}^Mf(mx)&=& \sum_{n=1}^\infty\m(n)\sum_{m=1}^Mg(mnx)\,=\, \sum_{n=1}^\infty\m(n)\sum_{m=1\atop nm\le M}^M
g(mnx) + R_M(x)
\end{eqnarray*}
where $R_M(x)={\displaystyle \sum_{n=1}^\infty}\m(n){\displaystyle\sum_{m=1\atop nm> M}^M}
g(mnx)$ verifies
\begin{eqnarray*}|R_M(x)|\,\le \,\sum_{n=1}^\infty{\displaystyle\sum_{m=1\atop nm> M}^M}
|g(mnx)|\,\le \, \sum_{  \nu> M}|g(\nu x)|\big(\sum_{ n\ge 1\atop n|\nu}1\big)\,= \, \sum_{  \nu> M}|g(\nu x)|d(\nu),
\end{eqnarray*}
and    thus is small if $M$ is large. Further, 
\begin{eqnarray*}\sum_{n=1}^\infty\m(n)\sum_{  nm\le M}g(mnx) \,=\, \sum_{  \nu \le M}g(\nu x)\Big( \sum_{n|\nu\atop n\ge \frac{\nu}{M}} \m(n) \Big)\,=\, \sum_{  \nu \le M}g(\nu x)\d(\nu)\,=\, g(x),
\end{eqnarray*}
by \eqref{inversemu}.

\vskip 2 pt  Note that \eqref{HWcondition}  is clearly satisfied if $g$ is finitely supported on integers. Next choose $f:[0,1]\to \R$ so that
$$ c_f(\ell)= \sum_{m=1}^\infty \m\big(m)D(m\ell).$$
where  \begin{eqnarray*} D(d)=\begin{cases} D_j\quad & \quad \hbox{if $d=d_j$,  $1\le j\le k$,}
\cr 0 \quad & \quad \hbox{otherwise.}
\end{cases}
\end{eqnarray*}
Then 
$$D(d)= \sum_{m=1}^\infty c_f(md)=D_f(d) .$$ 

Now, \begin{eqnarray*}F_n(f)&=&\sum_{d\le n}D_{f}(d)\,M\big(\frac{n}{d}\big)
\,= \, \sum_{1\le j\le k}D_{f}(d_j) \,M\Big(\frac{Y_{k+1}}{d_j}\Big)
\cr &=&\sum_{1\le j\le k}D_{f}(d_j)\,M(m_j)
\,\ge \,C\,\sum_{1\le j\le k}D_{f}(d_j) \sqrt {m_j}
\cr &\ge &C\,\sqrt n \sum_{1\le j\le k} \frac{D_f(d_j)}{\sqrt {d_j}}\,=\, C\,\sqrt n \sum_{d\le n} \frac{|D_f(d)|}{\sqrt {d}}.
\end{eqnarray*}

\end{proof}

\begin{remark} Assume that $g$ is positive everywhere. Let $x$ be  a  real such that for every $c>0$, (weaker conditions are available)
\begin{equation}\label{slicereg}  \g_1(c)=\sup_{   n\ge 1   }\sup_{ m\ge cn  }{g(m x) \over
g(n x)}  <\infty.
\end{equation}
Then as a special case of  Corollary 1.2 in \cite{W2}, 
condition \eqref{HWcondition} holds as soon as 
\begin{equation}\label{conv.} \sum_{n\ge 1}g(n x) \log n<\infty.
\end{equation}
In other words, it suffices that \eqref{HWcondition} holds with the mean value of $d(n)$ in   place of $d(n)$, since $\sum_{k\le n}  d(k) \sim n\log n$.\end{remark}
\begin{proof}[{\bf Proof of Theorem \ref{t1a}.}]  By Lemma \ref{l2},
\begin{equation*}F_{n,\s}(h)\ =\ \sum_{d\le n}\,\frac{R_{h_\s}(d)}{d^{2\s-1}}\,\Big( \sum_{ \l \le \frac{n}{d}}\frac{\m(\l)} {\l^{2\s}}\Big), 
\end{equation*}  
for any $\s\in \C$. Let $0\le \s <1$. We prove that
\begin{equation}\label{riemest}
\big|R_{h_\s}(\ell)- \int_{0}^{1} h_\s(t)\dd t\big| \le \frac{C_\s}{\ell^{1-\s}}.
\end{equation}
As   
$h_\s'(x)=(-2\pi \sin 2\pi x)x^{-\s}+\cos2\pi x(-\s x^{-1-\s})$, we have
  $|h_\s'(x)|\le (2\pi+\s)x^{-1-\s}$, for any $x\in]0,1]$. 
Let   $1\le k \le \ell$ and  $t\in ] \frac{k-1}{\ell}, \frac{k}{\ell}]$. Then, 
$$\big|h_\s(\frac{k}{\ell})-h_\s(t)\big|\ \le\ \frac{1}{\ell}|h_\s'(\xi)|\ \le\ \frac{1}{\ell} \  \frac{2\pi+\s}{\xi^{1+\s}},$$
for some $\xi\in] \frac{k-1}{\ell}, \frac{k}{\ell}[$.
Thus 
\begin{equation} \label{hcom}    \int_{ \frac{k-1}{\ell}}^{ \frac{k}{\ell}}  \big|h_\s(\frac{k}{\ell})-h_\s(t)\big|\dd t\,\le \,    \frac{C_\s}{\ell^2\,(\frac{k}{\ell})^{1+\s}} \ =\   \frac{C_\s}{\ell^{1-\s}}\, \frac{1}{k^{1+\s}}  ,\qq C_\s= 2^\s(2\pi+\s)
\end{equation}
 if $k>1$, and if $k=1$, 
$\int_0^{\frac{1}{\ell}} |h_\s(t)|\dd t\le \frac{1}{(1-\s)\ell^{1-\s}}$, and 
$\int_0^{\frac{1}{\ell}} |h_\s(\frac{1}{\ell})|\dd t\le \frac{1}{\ell^{1-\s}}$. Thus
\begin{eqnarray} \label{hsum} \sum_{k=1}^\ell   \int_{ \frac{k-1}{\ell}}^{ \frac{k}{\ell}}  \big|h_\s(\frac{k}{\ell})-h_\s(t)\big|\dd t&\le & \frac{C_\s}{\ell^{1-\s}}\,\sum_{k=1}^\ell   \frac{1}{k^{1+\s}}\ \le \ \frac{C_\s}{\ell^{1-\s}}  .
\end{eqnarray}

Now  
\begin{eqnarray} \label{hcomp} 
 \sum_{k=1}^\ell  \frac{1}{ k^{\s } }\, h\Big(  \frac{k}{\ell}\Big)-\ell^{1-\s} \int_{0}^{1} h_\s(t)\dd t&=&     \ell^{1-\s}\sum_{k=1}^\ell \int_{ \frac{k-1}{\ell}}^{ \frac{k}{\ell}}\Big[  \Big(\frac{k}{\ell }\Big)^{-\s}h\Big(  \frac{k}{\ell}\Big)-h_\s(t)\Big]\dd t
 \cr &=&    \ell^{1-\s}\sum_{k=1}^\ell \int_{ \frac{k-1}{\ell}}^{ \frac{k}{\ell}}\big[  h_\s\big(\frac{k}{\ell }\big)-h_\s(t)\big]\dd t
  .
\end{eqnarray}
Thus \begin{eqnarray*}  
\Big|\sum_{k=1}^\ell  \frac{1}{ k^{\s } }\, h\Big(  \frac{k}{\ell}\Big)- \ell^{1-\s}\int_{0}^{1} h_\s(t)\dd t\Big|
  \,\le\,  \ell^{1-\s}  \sum_{k=1}^\ell \int_{ \frac{k-1}{\ell}}^{ \frac{k}{\ell}}\big|  h_\s\big(\frac{k}{\ell }\big)-h_\s(t)\big|\dd t
  \,\le\,  C_\s  .
\end{eqnarray*}
Therefore $\big|R_{h_\s}(\ell)- \int_{0}^{1} h_\s(t)\dd t\big| \le {C_\s}/{\ell^{1-\s}}
$ as claimed.

\vskip 2 pt
We deduce that
  \begin{eqnarray} \label{happrox} F_{n,\s}(h)&=& \sum_{d\le n}\,\frac{R_{h_\s}(d)}{d^{2\s-1}}\,\Big( \sum_{ \l \le \frac{n}{d}}\frac{\m(\l)} {\l^{2\s}}\Big)
\cr &=&\Big(\int_{0}^{1} h_\s(t)\dd t\Big)\sum_{d\le n}\,\frac{1}{d^{2\s-1}}\,\Big( \sum_{ \l \le \frac{n}{d}}\frac{\m(\l)} {\l^{2\s}}\Big)
  + \mathcal O_\s\Big(\sum_{d\le n}\,\frac{1}{d^{\s}}\,\Big| \sum_{ \l \le \frac{n}{d}}\frac{\m(\l)}{\l^{2\s}}\Big|\Big).
\end{eqnarray} 
Let $1/2<\s <1$. By Lemma \ref{i-ii},
\begin{eqnarray*}
  \sum_{ d\le n}\frac{1}{d^{2\s-1}} \big(\sum_{\l\le \frac{n}{d}} \frac{\m(\l)}{\l^{2\s}}\big)\ =\   \frac{n^{2(1-\s)}} {2(1-\s)\zeta(2)}+ \frac{\zeta(2\s-1)}{\zeta(2\s)}+ \mathcal O \big(n^{1-2\s}\log n \big).
   \end{eqnarray*}
   Consequently,
\begin{eqnarray*}
 F_{n,\s}(h)&=&\Big(\int_{0}^{1} h_\s(t)\dd t\Big)\Big(\frac{n^{2(1-\s)}} {2(1-\s)\zeta(2)}+ \frac{\zeta(2\s-1)}{\zeta(2\s)}+ \mathcal O \big(n^{1-2\s}\log n \big)\Big)
 + \mathcal O_\s\big(n^{1-\s}\big)
 \cr &=& \frac{\int_{0}^{1} h_\s(t)\dd t } {2(1-\s)\zeta(2)}\,n^{2(1-\s)}+  \mathcal O_\s\big(n^{1-\s}\big).
\end{eqnarray*}
  
Now let $\s= 1/2$. By Lemma \ref{l2}, next   Lemma \ref{i-ii} and \eqref{A3.13},
\begin{eqnarray*}F_{n,1/2}(h)&=& \sum_{d\le n}\,R_{h_{1/2}}(d)\,\Big( \sum_{ \l \le \frac{n}{d}}\frac{\m(\l)} {\l}\Big), 
\cr &=&  \Big(\sum_{d\le n} \sum_{ \l \le \frac{n}{d}}\frac{\m(\l)} {\l}\Big)\int_{0}^{1} h_{1/2}(t)\dd t+\mathcal O\Big(\sum_{d\le n} \frac{1}{\sqrt d}\Big|\sum_{ \l \le \frac{n}{d}}\frac{\m(\l)} {\l}\Big|\Big)
\cr &=&  \Big(\frac{n}{\zeta(2)} + \mathcal O(\log n) \Big)\int_{0}^{1} h_{1/2}(t)\dd t+\mathcal O(\sqrt n).
 \cr &=& \Big(\int_{0}^{1} h_{1/2}(t)\dd t\Big)\frac{ n } {\zeta(2)}\, +  \mathcal O \big(\sqrt n\, \big).\end{eqnarray*}

Finally let $0<\s< 1/2$. By Lemma \ref{i-ii},
\begin{eqnarray*}
  \sum_{ d\le n}\frac{1}{d^{2\s-1}} \big(\sum_{\l\le \frac{n}{d}} \frac{\m(\l)}{\l^{2\s}}\big)&=& \frac{n^{2(1-\s)}} {2(1-\s)\zeta(2)}+ \mathcal O \big(n^{1-2\s}\log n \big),
 \end{eqnarray*}
and by \eqref{mobiusweighted}, $\sum_{ \l \le x}\frac{\m(\l)} {\l^{2\s}} =\mathcal O_\s (x^{1-2\s})$. Thus
\begin{eqnarray*}\sum_{ d\le n}\frac{1}{d^{\s }} \big(\sum_{\l\le \frac{n}{d}} \frac{\m(\l)}{\l^{2\s}}\big)\ =\   \mathcal O \Big(\sum_{ d\le n}\frac{1}{d^{\s }}\big(\frac{n}{d}\big)^{1-2\s}\Big)
\ =\    \mathcal O \big(n^{1-\s}\big).
   \end{eqnarray*}
By reporting these estimates in \eqref{happrox}, we get
 \begin{eqnarray*} F_{n,\s}(h) &=&\Big(\int_{0}^{1} h_\s(t)\dd t\Big)\Big(\frac{n^{2(1-\s)}} {2(1-\s)\zeta(2)}+ \mathcal O \big(n^{1-2\s}\log n \big)\Big) 
   + \mathcal O \big(n^{1-\s}\big)
   \cr &=&\Big(\int_{0}^{1} h_\s(t)\dd t\Big)\Big(\frac{n^{2(1-\s)}} {2(1-\s)\zeta(2)}\Big) 
   + \mathcal O \big(n^{1-\s}\big).
\end{eqnarray*} 

\end{proof}

\section{\bf Quadratic Riemann  sums.}\label{s3}
 We prove the following theorems.
 \begin{theorem}\label{p1}  Let $f:[0,1]\to \R$. Let $\frac{1}{2}<\s<1$  and assume that $f_\s\in \mathcal C$. 
 Further assume that the series 
 $$\sum_{\nu \in \Z}\Re\big(c_{f_\s}(\nu)\big),\qq \sum_{\nu\in \Z_* } \Re\big(c_{f_\s}(\nu)\big)
 \,\s_{1-2\s }(\nu)$$ are convergent, where $c_{f_\s}(\ell)$, $\ell \in \Z$ are the Fourier coefficients of ${f_\s}$, and $\s_{1-2\s }(\nu)$ denotes the sum of the $(1-2\s)$-th powers of the divisors of $\nu$. Then,
\begin{eqnarray*}S_{n,\s}(f)  &=& \big( \int_0^1 f_\s(x)\dd x\big)\Big(\frac{n^{2(1-\s)}}{2(1-\s)} + \zeta(2\s-1)-1\Big)
\cr & & \quad +\sum_{\nu\in \Z_* } \Re\big(c_{f_\s}(\nu)\big)\big( \s_{1-2\s }(\nu)- 1\big)+ o(1), 
\end{eqnarray*}
as $n\to \infty$.
\end{theorem}
 Theorem \ref{p1} will be deduced from the following preliminary result.
 
\begin{theorem}\label{p2} Let $\s$ be a real number. Let $f:[0,1]\to \R$  and assume that $f_\s\in \mathcal C$. 
 
 \vskip 2 pt {\rm (1)} We have,
 \begin{eqnarray*}  
 S_{n,\s}(f)&=& \big( \int_0^1 f_\s(x)\dd x\big)\Big( \sum_{1\le \ell \le n}(\ell -1)\ell^{-2\s }\Big)+ \sum_{\nu\in \Z_*}   c_{f_\s}(\ell)\Big(\sum_{1\le \ell \le n}\ell^{-2\s } \e_\nu (\ell)\Big),
\end{eqnarray*} 
 for each positive   $n$,    where  $ \e_\ell (d)$ is defined in \eqref{eps}.
\vskip 3 pt {\rm (2)}  Assume that the series $\sum_{\nu \in \Z}c_{f_\s}(\nu) $ is convergent. Then,
\begin{eqnarray*} S_{n,\s}(f)
 &=& \big( \int_0^1 f_\s(x)\dd x\big)\Big( \sum_{1\le \ell \le n}\ell^{1-2\s }\Big)+ \sum_{\nu\in \Z_* }  c_{f_\s}(\nu) \Big( \sum_{1\le \ell \le n\atop \ell|\nu}\ell^{1-2\s }\Big)
 \cr & & \quad - \sum_{\nu\in \Z} c_{f_\s}(\nu)\Big( \sum_{1\le \ell \le n}\ell^{-2\s }\Big),
\end{eqnarray*}
for each positive  $n$.\end{theorem}
  Without  assuming  $f_\s\in \mathcal C$, we have the following basic result.
\begin{theorem}\label{p3}   Assume that  
 $R_{f_\s}(\ell)$
converge to a finite limit  $I(f_\s)$, as $\ell\to \infty$.  Then
\begin{eqnarray*}
\lim_{n\to \infty} \frac{S_{n,\s}(f)}{\sum_{1\le \ell \le n}\ell^{1-2\s }} \ =\   I(f_\s).
\end{eqnarray*}
\end{theorem}

Consider   $g(x)$  defined in \eqref{g},
\begin{eqnarray*}    g(x)\ =\  \frac{\sin \log x}{\log x} \qq 0< x< 1,\qq g(1) =g(0)=0.
\end{eqnarray*}
  \begin{theorem}\label{t1}If $1/2< \s< 1$, then
\begin{eqnarray*}
\lim_{n\to \infty}\ \frac{ S_{n,\s}(g)}{n^{2(1-\s)}} \ =\  \frac{1}{2( 1-\s)}\, \arctan \big( \frac{1}{1-\s}\big).
\end{eqnarray*}
Further,
\begin{eqnarray*}
   \frac{ S_{n,\s}(g)}{n^{2(1-\s)}}\ = \ \frac{1}{2(1-\s)}\,  \arctan \big( \frac{1}{1-\s}\big)\, +\, \mathcal O_\s \Big( \frac{1}{n^{2(1-\s)/3}}\Big)\, .\end{eqnarray*}   \end{theorem}
\begin{remark} \label{sns.trivial bound} The asymptotic size's order of $ S_{n,\s}(g)$ is thus the one given by the trivial bound $S_{n,\s}(g)=\mathcal O_\s( n^{2(1-\s)})$.  
\end{remark}
As a corollary we get
\begin{corollary}\label{c1} Let $1/2< \s< 1$.  Then, 
\begin{eqnarray*}
 \frac{1}{n^{2(1-\s)}} \int_0^{1}\Big|\sum_{k=1}^n\frac{1}{k^{\s+it}}\Big|^2 \dd t\ =\  \frac{1}{1-\s}\,\arctan \big( \frac{1}{1-\s}\big) \, +\, \mathcal O_\s \Big( \frac{1}{n^{2(1-\s)/3}}\Big)\,,
\end{eqnarray*}
as $n$ tends to infinity.
\end{corollary}

\begin{remark} It will be clear from the proofs given, that the previous Theorems  extend with no difficulty to the modified Riemann quadratic sums
\begin{eqnarray}\label{S_n(f)modified} S_{n,\s,\a}(f) \ =\  \sum_{1\le k\le \ell \le n^\a}\frac{1}{(k\ell)^{\s }}\, f\Big( \frac{k}{\ell}\Big) 
\end{eqnarray}
where $0<\a <1$.
 \end{remark}


\subsection{Proof of Theorem \ref{p3}.}  
We note that
 \begin{eqnarray}\label{e3}
 S_{n,\s}(f) &=& \sum_{1\le k\le \ell \le n}\frac{1}{\ell^{2\s }}\, f_\s\Big( \frac{k}{\ell}\Big)
\ =\ \sum_{1\le \ell \le n}\ell^{1-2\s }\, R_{f_\s}(\ell),
\end{eqnarray}
where $R_{f_\s}(\ell)=D_{f_\s}(\ell)/\ell$. 
Thus 
\begin{eqnarray}\label{e31}
 \frac{S_{n,\s}(f) }{\sum_{1\le \ell \le n}\ell^{1-2\s }}&=&\sum_{1\le \ell \le n}\, v_{n,\ell}\, R_{f_\s}(\ell),\end{eqnarray}
where   \begin{eqnarray}\label{T0}
 v_{n,\ell}= \begin{cases} \ell^{1-2\s }\big\slash\big(\sum_{1\le \ell \le n}\ell^{1-2\s }\big)\quad &\hbox{if $1\le \ell \le n$},\cr 
 0 \quad &\hbox{otherwise}.
\end{cases} \end{eqnarray}
This reduces the problem to a matrix summation question. 
 In the next  lemma, we  just add  to well-known Toeplitz's criterion a rate of convergence. 
\begin{lemma}\label{[T]}
  Let  $\{ \t_{n,\ell},n\ge 1, 1\le \ell\le n\}$ be a triangular array of complex numbers verifying the following conditions
   \begin{eqnarray}\label{T}
   \begin{cases} {\rm (i)}\ \ \lim_{n\to \infty} \t_{n,\ell}\ =\ 0,\qquad \hbox{for each $\ell $},\cr 
{\rm (ii)}\ \ \sup_{n\ge 1} \sum_{ 1\le \ell \le n} |\t_{n,\ell} |\ =\ M<\infty .
\end{cases} \end{eqnarray} Let $\{x_{\ell}, \ell\ge 1\}$ be a bounded sequence of reals and set $T_n = \sum_{\ell=1}^n \t_{n,\ell} x_\ell$. Then for any $n\ge 1$ and any $D\ge 1$,
\begin{eqnarray}\label{bound}
\big|T_n\big| &\le & \big( \sup_{\ell \ge 1} |x_\ell|\big)\sum_{\ell \le D} |\t_{n,\ell}| + M\,\sup_{\ell >D}|x_\ell|.
\end{eqnarray}
In particular, if $\lim_{\ell \to \infty}x_\ell=0$, then $\lim_{n \to \infty}T_n=0$.\end{lemma}
\begin{proof} Immediate since
\begin{eqnarray*}
\big|T_n\big| &\le &  \sum_{\ell \le D} |\t_{n,\ell}||x_\ell|  +  \Big(\sup_{\ell >D}|x_\ell|\Big)\sum_{D<\ell \le n} |\t_{n,\ell}|
\cr &\le &\Big( \sup_{\ell \ge 1} |x_\ell|\Big)\sum_{\ell \le D} |\t_{n,\ell}| +  \Big(\sup_{\ell >D}|x_\ell|\Big)M.
\end{eqnarray*}
If $\lim_{\ell \to \infty}x_\ell=0$, given any positive real $\e$ and fixing $D=D(\e)$ sufficiently large so that $\sup_{\ell >D}|x_\ell|\le \e$, we have for any $n\ge 1$,
\begin{eqnarray*}
\big|T_n\big| &\le &   \big(\sup_{\ell \ge 1} |x_\ell| \big)\sum_{\ell \le D} |\t_{n,\ell}|  +  \e\, M
.
\end{eqnarray*}
 Whence  $\limsup_{n\to \infty} \big|T_n\big|\le  \e M
$, by \eqref{T}-(i). As $\e$ can be arbitrary small, this achieves the proof.\end{proof}
\vskip 7 pt
The triangular array  \eqref{T0} obviously verifies the  conditions \eqref{T}.
Assume that the limit  $I_\s(f)=\lim_{\ell\to \infty} R_{f_\s}(\ell) $ 
exists. By  \eqref{e31}, 
\begin{eqnarray}\label{e32}
 \frac{S_{n,\s}(f) }{\sum_{1\le \ell \le n}\ell^{1-2\s }}&=&\sum_{1\le \ell \le n}\, v_{n,\ell} R_{f_\s}(\ell) \cr&=& I_\s(f) + \sum_{1\le \ell \le n}\, v_{n,\ell}\big( R_{f_\s}(\ell) -I_\s(f)\big).
 \cr 
 \end{eqnarray}
  Lemma \ref{[T]} applied to $x_\ell =R_{f_\s}(\ell)-I(f_\s)$, $\ell \ge 1$,   together with  \eqref{e32} imply that 
\begin{eqnarray}\label{e32a}
 \frac{S_{n,\s}(f) }{\sum_{1\le \ell \le n}\ell^{1-2\s }}&=&I(f_\s) + o(1),\qq\qq n\to \infty,\end{eqnarray}
which achieves the proof.


\subsection{Proof of Theorem \ref{t1}.}   
 We  first prove  the convergence of $R_{g_\s}(\ell)$  to  $\int_0^1 g_\s(t) \dd t$, for any $\s\in [0,1[$, and provide a speed of convergence.

\begin{lemma}\label{l2R} We have 
\begin{eqnarray*}   \big|R_{g_\s}(\ell)- \int_{0}^{1} g_\s\dd t\big|&\le & \frac{C_\s}{\ell^{1-\s}}
.\end{eqnarray*}
Further\begin{eqnarray*}\label{} 
\int_{0}^{1} g_\s(t)\dd t &=&
\arctan \big( \frac{1}{1-\s}\big)
.
\end{eqnarray*}\end{lemma}

We need a lemma.

 \begin{lemma}\label{l2a} 
   We have 
 \begin{eqnarray*}   \sum_{1\le k\le \ell}  \int_{ \frac{k-1}{\ell}}^{ \frac{k}{\ell}}  \big|g_\s(\frac{k}{\ell})-g_\s(t)\big|\dd t&\le &    \frac{C_\s}{\ell^{1-\s}}  .
\end{eqnarray*}
\end{lemma}
\begin{proof}
 As $g'(x)=  {(\log x)\cos (\log x) - \sin(\log x) \over x(\log x)^2}$, we have
$$g_\s'(x)= \frac{g'(x)-\s x^{-1}g(x)}{x^\s}=\frac{1}{x^{1+\s}(\log x)}\Big\{\cos (\log x)-\frac{\sin(\log x)}{\log x} -\s \sin(\log x)\Big\}.$$
Thus  $|g_\s'(x)|\le  \frac{2+\s}{x^{1+\s}|\log x|}$, for any $x\in]0,1]$. 

Let  $\ell$ and $k$ be two arbitrary integers such that $1\le k \le \ell$, and let $t\in ] \frac{k-1}{\ell}, \frac{k}{\ell}]$. We have
$$\big|g_\s(\frac{k}{\ell})-g_\s(t)\big|\ \le\ \frac{1}{\ell}|g_\s'(\xi)|\ \le\ \frac{1}{\ell} \  \frac{2+\s}{\xi^{1+\s}|\log \xi|},$$
for some $\xi\in] \frac{k-1}{\ell}, \frac{k}{\ell}[$.
Thus if $k>1$,
\begin{eqnarray} \label{com}    \int_{ \frac{k-1}{\ell}}^{ \frac{k}{\ell}}  \big|g_\s(\frac{k}{\ell})-g_\s(t)\big|\dd t&\le &    \frac{C_\s}{\ell^2\,(\frac{k}{\ell})^{1+\s}\log\frac{\ell}{k}} \ =\   \frac{C_\s}{\ell^{1-\s}}\, \frac{1}{k^{1+\s}\log\frac{\ell}{k}}  .
\end{eqnarray}
If $k=1$,  then 
$\int_0^{\frac{1}{\ell}} |g_\s(t)|\dd t\le \int_0^{\frac{1}{\ell}} \frac{1}{t^\s \log \frac{1}{t}}\dd t\le \frac{1}{(1-\s)\ell^{1-\s}\log \ell}$, and $\int_0^{\frac{1}{\ell}} |g_\s(\frac{1}{\ell})|\dd t\le \frac{1}{\ell}\frac{\ell^\s}{\log \ell} =\frac{1}{\ell^{1-\s}\log \ell}$. So that the bound in \eqref{com} remains valid for $k=1$ either.
Thus
\begin{eqnarray} \label{sum} \sum_{k=1}^\ell   \int_{ \frac{k-1}{\ell}}^{ \frac{k}{\ell}}  \big|g_\s(\frac{k}{\ell})-g_\s(t)\big|\dd t&\le & \frac{C_\s}{\ell^{1-\s}}\,\sum_{k=1}^\ell   \frac{1}{k^{1+\s}\log\frac{\ell}{k}}   .
\end{eqnarray}
Now,
\begin{eqnarray*}      \sum_{1<k\le \ell/2}  \frac{1}{k^{1+\s}\log\frac{\ell}{k}} 
 &\le & C\,   \sum_{1<k\le \ell/2}  \frac{1}{k^{1+\s}}\ \le \   C_\s  .
 \end{eqnarray*}
If $\ell/2 < k\le \ell$, we write as in the proof of  lemma 7.2 in \cite{Ti},  $k=\ell-r$ where $1\le r< \ell/2$.  Then $\log \frac{\ell}{k}=-\log
\big(1-\frac{r}{\ell})>\frac{r}{\ell}$. 
Whence also,
\begin{eqnarray*} \sum_{\ell/2 < k\le \ell}  \frac{1}{k^{1+\s}\log\frac{\ell}{k}} 
&\le &    C_\s\sum_{1\le r< \ell/2}  \frac{1}{\ell^{1+\s}\,(\frac{r}{\ell})}\ \le \ C_\s\,  \frac{\log \ell}{\ell^{\s}}\, .
\end{eqnarray*}
By combining we  thus deduce
\begin{eqnarray}  \label{sigma} \sum_{1<k\le \ell}  \frac{1}{k^{1+\s}\log\frac{\ell}{k}} 
&\le &    C_\s .
\end{eqnarray}
By inserting \eqref{sigma} in \eqref{sum}, we therefore obtain
 \begin{eqnarray*}   \sum_{1<k\le \ell}  \int_{ \frac{k-1}{\ell}}^{ \frac{k}{\ell}}  \big|g_\s(\frac{k}{\ell})-g_\s(t)\big|\dd t&\le &    \frac{C_\s}{\ell^{1-\s}}\,  .
\end{eqnarray*}
\end{proof}
 
\begin{proof}[Proof of Lemma \ref{l2R}]
As for arbitrary positive integers  $\ell$ and $k$,
   \begin{eqnarray} \label{elemf}   
  \frac{1}{ k^{\s } }&=&   \ell^{1-\s} \int_{ \frac{k-1}{\ell}}^{ \frac{k}{\ell}}  \Big(\frac{k}{\ell }\Big)^{-\s}\dd t,
\end{eqnarray}
we   have
\begin{eqnarray} \label{comp} 
 \sum_{k=1}^\ell  \frac{1}{ k^{\s } }\, g\Big(  \frac{k}{\ell}\Big)&=&   \ell^{1-\s}\sum_{k=1}^\ell \int_{ \frac{k-1}{\ell}}^{ \frac{k}{\ell}}  \Big(\frac{k}{\ell }\Big)^{-\s}g\Big(  \frac{k}{\ell}\Big)\dd t
\cr &=&    \ell^{1-\s}\sum_{k=1}^\ell \int_{ \frac{k-1}{\ell}}^{ \frac{k}{\ell}}\Big[  \Big(\frac{k}{\ell }\Big)^{-\s}g\Big(  \frac{k}{\ell}\Big)-\frac{g(t)}{t^\s}\Big]\dd t
  +\ell^{1-\s} \int_{0}^{1} \frac{g(t)}{t^\s}\dd t
 \cr &=&    \ell^{1-\s}\sum_{k=1}^\ell \int_{ \frac{k-1}{\ell}}^{ \frac{k}{\ell}}\big[  g_\s\big(\frac{k}{\ell }\big)-g_\s(t)\big]\dd t
  +\ell^{1-\s} \int_{0}^{1} g_\s(t)\dd t.
\end{eqnarray}
By Lemma \ref{l2a}, 
 \begin{eqnarray*}   \sum_{1<k\le \ell}  \int_{ \frac{k-1}{\ell}}^{ \frac{k}{\ell}}  \big|g_\s(\frac{k}{\ell})-g_\s(t)\big|\dd t&\le &    \frac{C_\s}{\ell^{1-\s}}  .
\end{eqnarray*}
Thus \begin{eqnarray*}  
\Big|\sum_{k=1}^\ell  \frac{1}{ k^{\s } }\, g\Big(  \frac{k}{\ell}\Big)- \ell^{1-\s}\int_{0}^{1} g_\s(t)\dd t\Big|
  &\le&  \ell^{1-\s}  \sum_{k=1}^\ell\Big| \int_{ \frac{k-1}{\ell}}^{ \frac{k}{\ell}}\big[  g_\s\big(\frac{k}{\ell }\big)-g_\s(t)\big]\dd t\Big|
 \cr&\le &  \ell^{1-\s}  \sum_{k=1}^\ell \int_{ \frac{k-1}{\ell}}^{ \frac{k}{\ell}}\big|  g_\s\big(\frac{k}{\ell }\big)-g_\s(t)\big|\dd t
  \cr&\le &  C_\s  .
\end{eqnarray*}
Dividing both sides by $\ell^{1-\s}$, we get 
 \begin{eqnarray*} \label{WeigRiem1}   \Big|R_{g_\s}(\ell)- \int_{0}^{1} g_\s\dd t\Big|&\le & \frac{C_\s}{\ell^{1-\s}}
.\end{eqnarray*}
Whence the first part of the Lemma.
 As (Dwight \cite[863.4]{D}),
\begin{eqnarray}\label{D863.4} 
\int_0^\infty \frac{\sin t}{t} e^{-xt} \dd t &=& \arctan \frac{1}{x}, \qq \qq x>0,
\end{eqnarray}
we get (here comes the restriction $\s<1$),
\begin{equation}\label{gs} 
\int_{0}^{1}g_\s(t)\dd t \,=\, \int_{0}^{1} \frac{\sin (\log t)}{t^\s(\log t)}\dd t \,=\, \int_{0}^{\infty} \frac{\sin u}{u}e^{-(1-\s)u}\dd u
\,=\,
\arctan \big( \frac{1}{1-\s}\big).
\end{equation}
\end{proof}
 
 \begin{proof}[Proof of Theorem \ref{t1}] 
We get  in view of Theorem \ref{p3} and Lemma \ref{l2R},
 \begin{eqnarray*}
\lim_{n\to \infty} \ \frac{ S_{n,\s}(g)}{ \sum_{1 \le \ell\le n}\ell^{1-2\s }} \ =\  \arctan \big( \frac{1}{1-\s}\big).
\end{eqnarray*}
\vskip 3 pt

Now by \eqref{bound}, next by the first part of Lemma \ref{l2R}, for $1/2\le \s<1$,
\begin{eqnarray*}
\big|\sum_{1\le \ell \le n}\, v_{n,\ell}\big( R_{g_\s}(\ell)-\int_{0}^{1} g_\s(t)\dd t\big)\big| &\le & C_\s \Big( \frac{\sum_{\ell \le D} \ell^{1-2\s }}{{\sum_{1\le \ell \le n}\ell^{1-2\s }}} \Big)
\cr & &\quad + \,\sup_{\ell >D}\big| R_{g_\s}(\ell)-\int_{0}^{1} g_\s(t)\dd t\big|
\cr & \le & C_\s\Big( \big(\frac{D}{n}\big)^{2(1-\s)}+ \frac{1}{D^{1-\s}}\Big).
\end{eqnarray*}
Choosing $D\sim n^{2/3}$ gives
\begin{eqnarray*}
\big|\sum_{1\le \ell \le n}\, v_{n,\ell}\big( R_{\ell,\s}(g)-\arctan \big( \frac{1}{1-\s}\big)\big)\big| &\le & \frac{ C_\s}{n^{2(1-\s)/3}}\ .
\end{eqnarray*}
 As by \eqref{e32}
\begin{eqnarray*}
 \frac{S_{n,\s}(g) }{\sum_{1\le \ell \le n}\ell^{1-2\s }} &=&\int_{0}^{1} g_\s(t)\dd t + \sum_{1\le \ell \le n}\, v_{n,\ell}\Big( R_{g_\s}(\ell)- \int_{0}^{1} g_\s(t)\dd t\Big),
 \end{eqnarray*}
 we get
\begin{eqnarray}\label{e323}
 \frac{S_{n,\s}(g) }{\sum_{1\le \ell \le n}\ell^{1-2\s }}&=&\arctan \big( \frac{1}{1-\s}\big) + \mathcal O_\s\Big(\frac{1}{n^{2(1-\s)/3}}\Big).\end{eqnarray}
By \eqref{re.zetasum}, $
\sum_{1 \le \ell\le n}\ell^{1-2\s }
= 
   \frac{n^{2(1-\s)}}{2(1-\s)} + \zeta(2\s-1) + \mathcal O\left(n^{1-2\s}\right)$,
 so that
 \begin{eqnarray}\label{e324}
  \frac{2(1-\s) S_{n,\s}(g)}{n^{2(1-\s)}}&=& \frac{ \sum_{1\le \ell \le n}\ell^{1-2\s }}{n^{2(1-\s)}/2(1-\s)}\ \frac{S_{n,\s}(g) }{\sum_{1\le \ell \le n}\ell^{1-2\s }}
  \cr& = &
   \Big(1+ \frac{C_\s}{n^{2(1-\s)}}\Big)\, \Big(\arctan \big( \frac{1}{1-\s}\big) + \mathcal O_\s\big(\frac{1}{n^{2(1-\s)/3}}\big)\Big)\,
 \cr & =&\arctan \big( \frac{1}{1-\s}\big) + \mathcal O_\s\Big(\frac{1}{n^{2(1-\s)/3}}\Big)\, .
\end{eqnarray}
We have obtained
\begin{eqnarray*}
   \frac{ S_{n,\s}(g)}{n^{2(1-\s)}}\ = \ \frac{1}{2(1-\s)}\,  \arctan \big( \frac{1}{1-\s}\big)\, +\, \mathcal O_\s \Big( \frac{1}{n^{2(1-\s)/3}}\Big)\, .\end{eqnarray*}
  This achieves the proof.\end{proof}

\begin{remark}\label{sinealogx} Let $a\ge 1$ and let $h(x)= \frac{\sin a \log x}{\log x}$. The same proof allows one to get
\begin{eqnarray}\label{sinealogx1}
   \frac{ S_{n,\s}(h)}{n^{2(1-\s)}}\ = \ \frac{1}{2(1-\s)}\,  \arctan \big( \frac{a}{1-\s}\big)\, +\, a\,\mathcal O_\s \Big( \frac{1}{n^{2(1-\s)/3}}\Big)\, .\end{eqnarray}
 Indeed, in this case 
 $$h_\s'(x)=\frac{1}{x^{1+\s}(\log x)}\Big\{a\cos (a\log x)-\frac{\sin(a\log x)}{\log x} -\s \sin(a\log x)\Big\}.$$
Thus  $|h_\s'(x)|\le  \frac{2a+\s}{x^{1+\s}|\log x|}$, for any $x\in]0,1]$. So it suffices to substitute the new constant $2a+\s$ to $2+\s$ everywhere in the proof.
\vskip 3 pt
In the critical case  $a=n^\b$, this is better than  the classical bound $S_{n,\s}(g)=\mathcal O_\s( n^{2(1-\s)}\log n)$ (cf. \cite[Eq.\,(7.2.1)]{Ti}), only if $\b<{2(1-\s)/3}$. As $\s>1/2$, ${2(1-\s)/3}<1/3$; so $\b$ must be less than $1/3$.  See also  Remark \ref{sns.trivial bound}.
  \end{remark}
 
\begin{problem}\label{pbsinealogx} \rm   How to improve the error term in \eqref{sinealogx1} when $a=n^\b$ with $\b\ge 1/3$?
\end{problem}

\subsection{Proof of Corollary \ref{c1}.}
 Follows from Theorem \ref{t1} since
 \begin{align}\label{e1}\int_0^{1}\Big|\sum_{k=1}^n\frac{1}{k^{\s+it}}\Big|^2 \dd t 
 &  = \sum_{k=1}^n\frac{1}{k^{2\s}}+  2\sum_{1\le k<   \ell \le n}\frac{1}{(k\ell)^{\s }}\ \frac{\sin \log \frac{k}{\ell}}{\log \frac{k}{\ell}}
 .
    \end{align}
\subsection{Proof of Theorem \ref{p2}.}

 Here again we use  
 Dini's test (Bary \cite[Vol.\,I]{B})
which implies  that the Fourier series of ${f_\s}$
$$\s_{f_\s}(x)=\sum_{\nu\in \Z} c_{f_\s}(\nu)e^{2{\rm i}\pi\nu x},$$
converges to ${f_\s}(x)$ 
 for any real $x$, $0<x<1$.
Thus
  \begin{eqnarray*}  
{f_\s}\big(\frac{k}{\ell}\big)\ =\ c_{f_\s}(0)+ \sum_{\nu\in \Z_*} c_{f_\s}(\nu)e^{2{\rm i}\pi\nu\frac{k}{\ell}} , \qq \quad k=1, 2, \ldots, \ell-1.
\end{eqnarray*}
Whence 
\begin{eqnarray*}  
D_{f_\s}(\ell)&=& (\ell -1)c_{f_\s}(0)+ \sum_{\nu\in \Z_*} c_{f_\s}(\nu)\Big(\sum_{k=1}^{\ell -1}e^{2{\rm i}\pi\nu \frac{k}{\ell}}\Big)
\cr &=& (\ell -1)c_{f_\s}(0)+ \sum_{\nu\in \Z_*} c_{f_\s}(\nu)\e_\ell(\nu).
\end{eqnarray*}
 Therefore \begin{eqnarray*}  
S_{n,\s}(f) &=& \sum_{1\le \ell \le n}\ell^{-2\s }\, D_{{f_\s}}(\ell)
\cr &=& c_{f_\s}(0)\Big( \sum_{1\le \ell \le n}(\ell -1)\ell^{-2\s }\Big)+ \sum_{\nu\in \Z_*} c_{f_\s}(\ell)\Big(\sum_{1\le \ell \le n}\ell^{-2\s } \e_\nu (\ell)\Big).
\end{eqnarray*}
Whence the first assertion. 
 Now if the series $\sum_{\ell \in \Z}c_{f_\s}(\ell)$ is convergent, we can write
 \begin{eqnarray*} D_{f_\s}(\ell)&=&   (\ell-1) \, c_{f_\s}(0)+ \sum_{\nu\in \Z_*} c_{f_\s}(\nu)\Big(\sum_{k=1}^{\ell }e^{2{\rm i}\pi\nu \frac{k}{\ell}}-1\Big) \cr &=&   \ell \, c_{f_\s}(0)+ \ell \sum_{\nu\in \Z_*\atop \ell |\nu } c_{f_\s}(\nu)-\sum_{\nu\in \Z} c_{f_\s}(\nu) \end{eqnarray*}
since $\sum_{k=1}^{\ell }e^{2{\rm i}\pi\nu \frac{k}{\ell}}$ equals to $\ell$ or $0$ according to $\ell|\nu$ or not. 
Thus 
\begin{eqnarray*}& & S_{n,\s}(f) \ =\ \sum_{1\le \ell \le n}\ell^{-2\s }\, D_{{f_\s}}(\ell)
\cr
 &=& c_{f_\s}(0)\Big( \sum_{1\le \ell \le n}\ell^{1-2\s }\Big)+ \sum_{\nu\in \Z_* } c_{f_\s}(\nu)\Big( \sum_{1\le \ell \le n\atop \ell|\nu}\ell^{1-2\s }\Big)
  - \sum_{\nu\in \Z} c_{f_\s}(\nu)\Big( \sum_{1\le \ell \le n}\ell^{-2\s }\Big).\end{eqnarray*}

  \subsection{Proof of Theorem \ref{p1}.} As $f$ is real-valued
\begin{eqnarray*}  S_{n,\s}(f)&=& c_{f_\s}(0)\Big( \sum_{1\le \ell \le n}\ell^{1-2\s }\Big)+ \sum_{\nu\in \Z_* }  \Re\big(c_{f_\s}(\nu)\big)\Big( \sum_{1\le \ell \le n\atop \ell|\nu}\ell^{1-2\s }\Big)
\cr & &\quad   - \sum_{\nu\in \Z} \Re\big(c_{f_\s}(\nu)\big)\Big( \sum_{1\le \ell \le n}\ell^{-2\s }\Big).
\end{eqnarray*}For each $N\ge 1$,
 \begin{eqnarray*} \lim_{n\to \infty} \sum_{-N\le \nu\le N\atop
\nu\neq 0 } \Re\big(c_{f_\s}(\nu)\big)\Big( \sum_{1\le \ell \le n\atop \ell|\nu}\ell^{1-2\s }\Big)& =&  \sum_{-N\le \nu\le N\atop
\nu\neq 0 }\Re\big(c_{f_\s}(\nu)\big)
 \s_{1-2\s }(\nu).
\end{eqnarray*}
From  the assumptions made, by using a standard approximation argument,  we deduce 
\begin{eqnarray*} \lim_{n\to \infty} \sum_{\nu\in \Z_* } c_{f_\s}(\nu)\Big( \sum_{1\le \ell \le n\atop \ell|\nu}\ell^{1-2\s }\Big)& =& \sum_{\nu\in \Z_* } c_{f_\s}(\nu)
 \s_{1-2\s }(\nu).
\end{eqnarray*} Using Theorem \ref{p2} and  estimate \eqref{re.zetasum} we get,
\begin{eqnarray*}S_{n,\s}(f)  &=& c_{f_\s}(0)\Big(\frac{n^{2(1-\s)}}{2(1-\s)} + C_\s\Big)+\sum_{\nu\in \Z_* } \Re\big(c_{f_\s}(\nu)\big)
 \s_{1-2\s }(\nu)- \sum_{\nu\in \Z} \Re\big(c_{f_\s}(\nu)\big)+ o(1), 
\end{eqnarray*}
as $n\to \infty$.

\vskip 5 pt

  

In the next Proposition we  provide with \eqref{varpi2} an $\ell^1$-type control of   the error term  
$$ \Theta_{n,\s}(f)= S_{n,\s}(f) - \big(\int_0^1 f_\s(t)\dd t\big)\sum_{1\le \ell \le n}\ell^{1-2\s }.$$
\begin{proposition}\label{varpi}Let $\bar f_\s= f_\s-\int_0^1 f_\s(t)\dd t$. Let $\varpi (n)$ be positive reals such that the series $\sum_{n\ge 1} 1/\varpi (n)$ converges, and let $\rho(d)= \sum_{n\ge d} 
 1/\varpi (n)$. 
 Assume that 
 \begin{equation}\label{varpi1} \sum_{ d \ge 1}d^{1-2\s }\rho(d)\, |R_{\bar f_\s}(d)|<\infty.
 \end{equation}
 Then
 \begin{equation}\label{varpi2}\sum_{n\ge 1} \frac{|\Theta_{n,\s}(f) |}{\varpi (n)}<\infty.
 \end{equation}
\end{proposition}
 \begin{proof}   By  linearity, 
 $R_{\bar f_\s}(\ell)=R_{ f_\s}(\ell)-\int_0^1 f_\s(t)\dd t$, next by \eqref{e3}, 
  \begin{eqnarray*}
 S_{n,\s}(f) &=& \sum_{1\le \ell \le n}\ell^{1-2\s }\, \Big(R_{\bar f_\s}(\ell)-\int_0^1 f_\s(t)\dd t+\int_0^1 f_\s(t)\dd t\Big)
 \cr &=& \Big(\int_0^1 f_\s(t)\dd t\Big)\Big(\sum_{1\le \ell \le n}\ell^{1-2\s }\Big)+ \sum_{1\le \ell \le n}\ell^{1-2\s }\,R_{\bar f_\s}(\ell).
\end{eqnarray*}
Thus
 \begin{eqnarray*}
\sum_{n\ge 1} \frac{|T_{n,\s}(f) |}{\varpi (n)}&=& \sum_{n\ge 1} \frac{\big|\sum_{1\le d\le n}d^{1-2\s }\,R_{\bar f_\s}(d)\big|}{\varpi (n)}\, \le
\, \sum_{n\ge 1} \frac{1}{\varpi (n)}\sum_{1\le d \le n}d^{1-2\s }\, |R_{\bar f_\s}(d)|
\cr &\le &
\sum_{ d \ge 1}d^{1-2\s }\, |R_{\bar f_\s}(d)|\sum_{n\ge d} \frac{1}{\varpi (n)}
\cr &\le &
\sum_{ d \ge 1}d^{1-2\s }\rho(d)\, |R_{\bar f_\s}(d)|\, <\, \infty,
\end{eqnarray*}
by assumption. 
 \end{proof}
\begin{remark}In the case of example \eqref{g}, it can be easily checked that Proposition \ref{varpi} however provides a weaker estimate than (see in particular \eqref{e323}) the one proven in Theorem \ref{t1}.
\end{remark}

\section{
Concluding remarks:  Local integrals of $\boldsymbol{\zeta(s)}$ and  amalgams.
}\label{s4}
We discuss some questions  related to the  previous sections, in particular to Remark \ref{sinealogx}, and to local integrals of the $\zeta$-function. 
We notably consider three interesting  related problems. Let $1\le a<b<\infty$, $a\ge cb$ for some positive $c$.  
We first note 
 that  for $\s>\frac12$,
\begin{align}\label{local.int.zeta}
\int_a^{b} |\zeta(\s +it ) |^2 \dd t  =-\zeta(2\s)(b-a)+ \,4S_{b,\s}(g(a,b)) +{\mathcal O}_{1/2,c}\,(
b^{1-2\s}),
\end{align}
where 
$$g(a,b)(x)= \frac{\sin(\frac12( b-a)\log x)\cos(\frac12( a+b)\log x)}{\log x}\qq\quad  0<x<1,$$
and $g(a,b)(1)=\frac12( b-a)$, $g(a,b)(0)=0$.
Indeed (as $\sin a -\sin b= 2\sin \frac12(a-b) \cos \frac12(a+b)$\,),
 \begin{align}\label{e1a}\int_a^{b}&\Big|\sum_{k=1}^b\frac{1}{k^{\s+it}}\Big|^2 \dd t 
  \  = (b-a)\sum_{k=1}^b\frac{1}{k^{2\s}}+  2\sum_{1\le k<   \ell \le b}\frac{1}{(k\ell)^{\s }}\ \frac{\sin b\log \frac{k}{\ell}-\sin  a\log \frac{k}{\ell}}{\log \frac{k}{\ell}}
 \cr  & =(b-a) \sum_{k=1}^b\frac{1}{k^{2\s}}+  2\sum_{1\le k<   \ell \le b}\frac{1}{(k\ell)^{\s }}\ \frac{\sin(\frac12( b-a)\log \frac{k}{\ell})\cos(\frac12( a+b)\log \frac{k}{\ell})}{\frac12 \log \frac{k}{\ell}}
   \cr 
   &=-(b-a)\sum_{k=1}^b\frac{1}{k^{2\s}}+ \,2S_{b,\s}(2g(a,b)). 
    \end{align}
By   the  classical approximation formula (\cite{Te}, Theorem 3.5), given $\s_0>0$, $0<\d<1$, we have uniformly for $\s\ge \s_0$, $x\ge 1$ , $0<|t| \le (1-\d)2\pi x$,
\begin{equation}\label{approx} \zeta(s) =\Sigma_x(s) -{x^{1-s}\over 1-s} +{\mathcal O}(
x^{-\s}), \qq \quad \Sigma_x(s)= \sum_{k\le x} {1\over k^s}.
\end{equation}
Thus
\begin{align*}
\int_a^{b} |\zeta(\s +it ) |^2 \dd t &=-(b-a)\sum_{k=1}^b\frac{1}{k^{2\s}}+ \,2S_{b,\s}(2g(a,b)) +{\mathcal O}_{1/2,c}\,(
b^{1-2\s})
\cr &=-\zeta(2\s)(b-a)+ \,4S_{b,\s}(g(a,b)) +{\mathcal O}(
b^{1-2\s}).
\end{align*}
\vskip 10 pt Now assume  $a=n$ and $b=n+1$. Then
\begin{equation*}\int_n^{n+1} |\zeta(\s +it ) |^2 \dd t\,=\, \hbox{{\it Quadratic Riemann sum} + {\it Error term}}\,.
\end{equation*}
More precisely,
\begin{align}\label{zeta.quadraticRiemann}
\int_n^{n+1} |\zeta(\s +it ) |^2 \dd t \ =\ -\zeta(2\s)+ \,4S_{n,\s}(g_{n}) +{\mathcal O}(
n^{1-2\s}),
\end{align}
where
\begin{eqnarray}\label{tildeg}g_n(x) =\frac{\cos((n+\frac12)\log x)\sin \frac12\log x}{\log x} \qq \hbox{if $0< x< 1$,}
\end{eqnarray}
  and $g_n(1) =\frac12$, $g_n(0)=0$. Let also  $g_{n,\s}(x)=g_n(x)/x^\s$,    $n\ge 1$. 
 \vskip 3 pt 
This stresses if necessary,  the importance of the quadratic sums $S_{n,\s}(g_{n}) $.  We could not  find  in the  literature asymptotic estimates for these sums when $n\to \infty$, even partial ones, which is a bit surprising.  
 The key quantities to be estimated are the Riemann sums $R_{g_{n,\s}}(d)$, $d\le n$. However  $g_{n,\s}(x)$    oscillates wildly  near $0_+$, with   peaks increasing to infinity with $n$,
and it seems illusory to directly estimate them. The Fourier coefficients can  however be computed.
We first 
  note by arguing as in \eqref{sfsigma} that
 \begin{equation}\label{tildeg1}\sum_{\nu\in \Z_*}|c_{g_{n,\s}}(\nu) |=\infty \ \  \hbox{and}\ \   \sum_{\nu\in \Z_*}|c_{g_{n,\s}}(\nu) |^q<\infty, \quad (q> \frac{1-\s}{\s},  \   n\ge 1).
  \end{equation}
 
 Next  Theorem \ref{p2} directly implies the following Proposition.  
 \begin{proposition}\label{prop.thp2}for $\frac{1}{2}<\s<1$,  
\begin{equation}\label{z1a}\int_n^{n+1} |\zeta(\s +it ) |^2 \dd t  =  -  \zeta(2\s)+ 4\sum_{k\in \Z_*} \Re\big( c_{g_{n,\s}}(k)\big)\Big(\sum_{1\le \ell \le n}\ell^{-2\s } \e_k (\ell)\Big)   + \mathcal O_\s\big(n^{1-2\s }\big),
\end{equation} 
 as $n\to \infty$, recalling that  $\e_k (\ell)$ is defined in \eqref{eps}.
Further  for any $k\in \Z*$, 
\begin{equation} \label{z1c}
\Re (c_{g_{n,\s}}(k) )\,=\,
 \sum_{m=0}^\infty \frac{(-1)^m(2 \pi k)^{2m}}{(2m)!}\Big(\frac{\arctan\frac{n+1}{2m+1-\s}-\arctan\frac{n}{2m+1-\s}}{2}\Big)
.\end{equation}
\end{proposition}
 \begin{proof}
By \eqref{zeta.quadraticRiemann},
\begin{align*}\label{zeta.quadraticRiemann}
\int_n^{n+1} |\zeta(\s +it ) |^2 \dd t \ =\ -\zeta(2\s)+ \,4S_{n,\s}(g_{n}) +{\mathcal O}(
n^{1-2\s}).
\end{align*}
Obviously $g_{n,\s}(x)$ is derivable   for any real $x$, $0<x<1$.  Thus by Theorem \ref{p2},
 \begin{eqnarray*}  
 S_{n,\s}(g_{n})\,=\, c_{g_{n,\s}}(0)\Big( \sum_{1\le \ell \le n}(\ell -1)\ell^{-2\s }\Big)+ \sum_{k\in \Z_*} \Re\big( c_{g_{n,\s}}(k)\big)\Big(\sum_{1\le \ell \le n}\ell^{-2\s } \e_k (\ell)\Big).
\end{eqnarray*}

Whence
\begin{align*}
\int_n^{n+1} |\zeta(\s +it ) |^2 \dd t &=-\zeta(2\s)+ \,4c_{g_{n,\s}}(0)\Big( \sum_{1\le \ell \le n}(\ell -1)\ell^{-2\s }\Big)+ \cr &\quad \,4\sum_{k\in \Z_*} \Re\big( c_{g_{n,\s}}(k)\big)\Big(\sum_{1\le \ell \le n}\ell^{-2\s } \e_k (\ell)\Big) +{\mathcal O}(
n^{1-2\s}).
\end{align*}
  
  Note that $c_{g_{n,\s}}(0)$ is  small, $c_{g_{n,\s}}(0)\sim \,{C_\s}/{n^2}$, since 
\begin{eqnarray*}
c_{g_{n,\s}}(0) &=&\int_0^1g_{n,\s}(t) \dd t  \,=\,
\frac{\arctan\frac{n+1}{1-\s}-\arctan\frac{n}{1-\s}}{2}
\sim \,\frac{C_\s}{n^2}.\end{eqnarray*}
Thus 
 \begin{eqnarray*}
c_{g_{n,\s}}(0)\Big( \sum_{1\le \ell \le n}(\ell -1)\ell^{-2\s }\Big)
=\mathcal O_\s\big(n^{-2\s }\big) .\end{eqnarray*}
Therefore 
\begin{eqnarray*}\int_n^{n+1} |\zeta(\s +it ) |^2 \dd t  &=& - \,\zeta(2\s)+ 4\sum_{k\in \Z_*} \Re\big( c_{g_{n,\s}}(k)\big)\Big(\sum_{1\le \ell \le n}\ell^{-2\s } \e_k (\ell)\Big)
\cr & &\qq + \mathcal O_\s\big(n^{1-2\s }\big),
\end{eqnarray*} 
as claimed. 
\vskip 2 pt Now we compute the Fourier coefficients of $g_{n,\s}$.  We have 
 \begin{eqnarray*}\label{} 
  c_{g_{n,\s}}(k)&=& \int_{0}^{1} g_{n,\s}(x)e^{-2{\rm i} \pi kx}\dd x
    \cr
&=&   \sum_{j=0}^\infty \frac{(-2{\rm i} \pi k)^j}{j!}\int_{0}^{1}  \frac{\sin (n+1)\log x-\sin  n\log x}{2x^\s\log x}\,x^j\,\dd x
\cr
 &=& \sum_{j=0}^\infty \frac{(-2{\rm i} \pi k)^j}{j!}\int_{0}^{\infty}\Big(e^{-(\frac{j+1-\s}{n+1})v}-e^{-(\frac{j+1-\s}{n})v}\Big)\frac{\sin v}{2 v}\dd v
\cr &=&  \sum_{j=0}^\infty \frac{(-2{\rm i} \pi k)^j}{j!}\Big(\frac{\arctan\frac{n+1}{j+1-\s}-\arctan\frac{n}{j+1-\s}}{2}\Big).
\end{eqnarray*}
Thus 
\begin{equation} \label{seriescoeff}
\Re (c_{g_{n,\s}}(k) )\,=\,
 \sum_{m=0}^\infty \frac{(-1)^m(2 \pi k)^{2m}}{(2m)!}\Big(\frac{\arctan\frac{n+1}{2m+1-\s}-\arctan\frac{n}{2m+1-\s}}{2}\Big)
.\end{equation}
\end{proof}
 
 \vskip 9 pt 
We conclude this paper with a remark on the $\mathcal S^2 $-Stepanov norm of $\zeta(s)$.
The Stepanov space $\mathcal S^2$ is   defined as the sub-space of functions $f$ of $L^2_{\rm loc}(\R)$     verifying the   following analogue of Bohr almost periodicity property: {\it For all $\e>0$, there exists $K_\e>0$ such that for any $x_0\in \R$, there exists $\tau\in [x_0, x_0+K_\e]$ such that  $\| f(.+\tau)-f(.)\|_{\mathcal S^2} \le \e$}. The Stepanov norm in  $\mathcal S^2 $ is equivalent to the \lq\lq amalgam\rq\rq\,norm  
  $$\|f\|_{\mathcal S^2}=\sup_{x\in \R}\Big( \int_x^{x+1} |f(t)|^2\dd t\Big)^{1/2}. $$
  The  almost everywhere convergence properties of   almost periodic Fourier series  in  $\mathcal S^2$ and of corresponding series of dilates were recently studied with Cuny in \cite{CW}; and a new form of Carleson's theorem for almost periodic Fourier series was proved. 
\vskip 2 pt A natural question arising from this study concerns the Riemann zeta function $\zeta(s )$, $s=\s +it$, and more precisely the evaluation of its Stepanov's norm, namely the supremum over all $n$ of the local integrals $\int_n^{n+1}|\zeta(\s +it ) |^2 \dd t$. 
It is clear that
\begin{lemma}\label{demi}
\begin{equation*} \big\|\zeta\big(\frac{1}{2} +i\cdot \big)\big\|_{\mathcal S^2}\ =\ \infty.\end{equation*} 
\end{lemma}
Indeed otherwise it would imply that
 $\int_0^T |\zeta({1\over 2} +it)|^2 \dd t \le C\,  \|\zeta\big(\frac{1}{2} +i\cdot \big)\|_{\mathcal S^2} T$,
which     contradicts the classical mean value estimate
$\int_0^T |\zeta({1\over 2} +it)|^2 \dd t \sim T\log T$,
(\cite{Ti} p.\,176).    
 Note that  
  \begin{equation}\label{omega.zeta.stepanov}\int_n^{n+1} \big|\zeta({1\over 2}  +it ) \big|^2 \dd t=\O \big(\log  n\big).
   \end{equation}
 
 It is interesting to observe that the above $\Omega$-result is in a sense optimal when  $t$ is
modelled by a Cauchy random walk. The behavior of the Riemann zeta-function on the critical line, along the Cauchy random walk, was studied by
Lifshits and Weber in \cite{LW}. Let $X_1,X_2,\ldots$ denote an
infinite sequence of independent Cauchy distributed random variables
(with characteristic function $\varphi(t) =e^{-|t|}$), 
and consider the partial sums
$S_n= X_1+ \ldots +X_n $. We shall prove the following 
precise result.

\begin{proposition}\label{propCauchy.increments} We have
\begin{eqnarray*}
\E  \Big| \zeta\Big({1\over 2} + iS_{n}\Big)\Big|^2\,=\,\mathcal O (\log n).\end{eqnarray*}
Moreover,\begin{eqnarray*}\E  \Big| \zeta\Big({1\over 2} + iS_{n+2}\Big)-\zeta\Big({1\over 2} + iS_{n}\Big)\Big|^2&=&2\,\log n +o(\sqrt{\log n}), \qquad n\to\infty.
\end{eqnarray*}
\end{proposition}
\begin{proof} Put
\begin{equation} \label{e118}
  \zeta_{1/2,n}=\zeta({1\over 2} + iS_n), \quad\qq n=1,2, \ldots 
\end{equation}
Put for any positive integer $n$
\begin{equation} \label{e119}
{\mathcal Z}_n = {\zeta}({1/  2}+iS_n) -\E  {\zeta}({1/ 2}+iS_n)
=\zeta_{1/2,n}-\E \zeta_{1/2,n}. 
\end{equation}
By Theorem 1.1 in Lifshits and Weber \cite{LW}, there exist explicit constants $C,C_0$ such that
\begin{eqnarray}\label{th.LW}  
& &{\rm (i)}\quad \E  \big|{\mathcal Z}_n \big|^2
\,=\,  \log n + C + o(1), \qquad n\to\infty,
\cr
& &{\rm (ii)}\quad \hbox{\rm For $m>n+1$, \qq}\big|\E\,{\mathcal Z}_n\overline{{\mathcal Z}_m}\big|
\,\le \,  C_0 \max \Big( {1\over n},{1\over 2^{ m-n  }}   \Big). \qq
\end{eqnarray}
Therefore, 
\begin{eqnarray}\label{e.diff}
  \E  \big|{\mathcal Z}_{n+2}-{\mathcal Z}_n \big|^2
&=& \E  \big|{\mathcal Z}_{n+2}\big|^2+\E  \big|{\mathcal Z}_{n}\big|^2- 2 \E{\mathcal Z}_n\overline{{\mathcal Z}_{n+2}}\
\cr &= &  \log (n+2)+ \log n +o(1)- 2 \E{\mathcal Z}_n\overline{{\mathcal Z}_{n+2}} \
\cr &= & 2\log n +\mathcal O(1).
\end{eqnarray}
By \cite[(3.3)]{LW},
\begin{equation} \label{e33}
\E  \zeta_{1/2,n}    =\zeta({1\over 2}+n)
-   {8n\over 4n^2-1}\,=\, 1+o(1).
\end{equation}
By combining with \eqref{th.LW}-(i) the first claim follows. As $|a-b|^2\le 2(a^2 +b^2)$, we also get
\begin{eqnarray*}\E  \big|(\zeta_{1/2,n+2}-\zeta_{1/2,n} )\big|^2&=& \E  \big|{\mathcal Z}_{n+2}-{\mathcal Z}_n +(\E \zeta_{1/2,n+2}-\E \zeta_{1/2,n} )\big|^2
\cr &=& \E  \big|{\mathcal Z}_{n+2}-{\mathcal Z}_n +o(1)\big|^2
\cr &=& 2\,\log n +o(\sqrt{\log n}),
\end{eqnarray*}
which proves the second claim.\end{proof}
For $\frac12<\s <1$,  the simple argument used in Lemma \ref{demi} no longer works since  by  Landau and Schnee's result, 
  $\int_0^T |\zeta(\s+it)|^2 \dd t\sim \zeta(2\s) T$  if $\s> \frac12$; in fact (\cite[p.13\&16]{Ma})
  $$\int_0^T |\zeta(\s+it)|^2 \dd t=\zeta(2\s)T +(2\pi)^{2\s-1}\frac{\zeta\big((2(1-\s)\big)}{2(1-\s)}T^{2(1-\s)} + \mathcal O\big( T^{\frac{2(1-\s)}{3}}(\log T)^{\frac{2}{9}}\big),$$
 if $\frac12<\s <1$.   
Thus  
$$\int_n^{n+1} \big|\zeta(\s  +it ) \big|^2 \dd t\ge A ,$$
for infinitely many $n$, with $A<\zeta(2\s)$, and the question arises whether $A$ can be arbitrary large, namely whether
\begin{problem}\label{prob2}  It is true that for $1/2<\s\le 1$
$$\limsup_{n\to \infty} \int_n^{n+1} \big|\zeta(\s  +it ) \big|^2 \dd t=\infty \,?$$  
\end{problem}
This   question is of a different kind from that one of   estimating the local integrals $\int_n^{n+1} \big|\zeta(\s  +it ) \big|^2 \dd t$;   both problems might rely on independent devices.
\vskip 3 pt
Maybe the answer to Problem \ref{prob2} is negative. In this direction, it seems that Proposition \ref{propCauchy.increments} can  be extended to the range of values $1/2<\s<1$ as follows:
\begin{eqnarray*}\E  \big| \zeta\big(\s + iS_{n+2}\big)-\zeta\big(\s+ iS_{n}\big)\big|^2&=&\mathcal O(1), \qquad n\to\infty.
\end{eqnarray*}
This is examined in a separate work \cite{W2a}. 
\vskip 3 pt

We also   note the following equivalent reformulation: for any $\s>0$, any $b>0$,
 \begin{equation}\label{step.norm.reformulation}k_b\, \|\zeta (\s+ i.)\|_{\mathcal S^2}\ \le\ \sup_{u\in \R}\ \int_\R  \frac{|\zeta (\s+ iv)|^2 }{b^2+ (v-u)^2}\dd v\ \le\ K_b\, \|\zeta (\s+ i.)\|_{\mathcal S^2}^2,
  \end{equation}
where $k_b, K_b$ are positive finite constants depending on $b$ only.  This follows from the Lemma below.

\begin{lemma}\label{l2loc} Let $b>0$. For any $f\in L^2_{{\rm loc}}(\R)$,
$$\|f \|_{\mathcal S^2}<\infty \quad \hbox{if and only if} \quad\sup_{u\in \R}\ \int_\R |f(t+u)|^2  \frac{\dd t}{b^2+ t^2}<\infty.$$
Moreover,
\begin{eqnarray*}k_b\, \|f\|_{\mathcal S^2}\ \le\ \sup_{u\in \R}\ \int_\R |f(t+u)|^2  \frac{\dd t}{b^2+ t^2}\ \le\ K_b\, \|f\|_{\mathcal S^2}^2, \end{eqnarray*}
where $k_b =\frac{1}{b^2+ 1}$ and $K_b=2\sum_{w\ge 0} \frac{1}{ b^2+ w^2} + \frac{1}{ b^2}$.
\end{lemma}

\begin{proof} On the one hand,
for any real $u$,
\begin{eqnarray*}\int_\R |f(t+u)|^2  \frac{\dd t}{b^2+ t^2}\ \ge \ \frac{1}{b^2+ 1} \int_0^1 |f(t+u)|^2  \dd t
\ = \  \frac{1}{b^2+ 1}\int_u^{u+1} |f(v)|^2   \dd v. \end{eqnarray*}
Taking supremum over $u$ in both sides yields
\begin{eqnarray*}\sup_{u\in \R}\ \int_\R |f(t+u)|^2  \frac{\dd t}{b^2+ t^2}&\ge &  \frac{1}{b^2+ 1}\, \|f\|_{\mathcal S^2}. \end{eqnarray*}
On the other hand,
 \begin{eqnarray*}\int_\R |f(t+u)|^2\frac{\dd t}{b^2+ t^2}&= &\sum_{x\in \Z}\int_{x}^{x+1} |f(v)|^2 \frac{\dd v}{b^2+ (v-u)^2}
 \cr&\le & \|f\|_{\mathcal S^2}^2\,\sum_{y\in \Z} \max_{y\le v\le y+1} \Big(\frac{1}{ b^2+ (v-\{u\})^2}\Big). \end{eqnarray*}
 


 
\end{proof}

Concerning the  integrals appearing in \eqref{step.norm.reformulation}, we recall   the classical formula \cite[(2.1.5)]{Ti},\begin{equation}\label{z1b}
\frac{\zeta(s)}{s} =- \int_0^\infty \{x\}x^{-1-s}\, \dd x=- \int_0^\infty \big\{\frac{1}{x}\big\}x^{s-1} \, \dd x,\qq\qq 0<\s<1 ,
\end{equation} 
where $\{x\}=x-[x]$ denotes the fractional part of $x$. Thus by Parseval equality for Mellin's transform, 
  \begin{eqnarray}\label{Parseval-zeta}
 \frac1{2\pi} \int_{-\infty}^{+\infty} \frac{\big|\zeta(\s +it)\big|^2}{\s^2+t^2} \dd t&=&\int_0^\infty \big\{\frac{1}{x}\big\}^2x^{2\s-1}\dd x
 \cr& =& -\frac{1}{\s} \, \Big(\frac{\zeta(2\s)}{2}+\frac{\zeta(2\s-1)}{2\s-1}\Big) ,
\end{eqnarray}
for $ 0<\s<1$, the   quantity in brackets being negative.  For the first equality, see \cite[Cor.\,1]{I}. The second equality was proved by Coffey  \cite[Prop.\,7]{C}, as well as   \cite[Cor.\,8]{C},
 \begin{eqnarray}\label{Parseval-zeta-half}
 \frac1{2\pi} \int_{-\infty}^{+\infty} \frac{\big|\zeta(\frac1{2} +it)\big|^2}{ \frac1{4}+t^2} \dd t&=&\log (2\pi) -\g,\end{eqnarray}
and we note   that  by Lemma \ref{demi} and  Corollary \ref{c1},
 \begin{eqnarray}\label{unif-zeta-half} \sup_{u\in \R}\ \ \int_{-\infty}^{+\infty}  \frac{|\zeta (\frac1{2}+ it)|^2 }{\frac1{4}+ (t-u)^2}\dd t\ =\ \infty.
  \end{eqnarray} 
\begin{remark} \label{Parseval-zeta-remark}
It is not clear to the author whether the integrals $$\int_{-\infty}^{+\infty}  \frac{|\zeta (\s+ iv)|^2 }{\s^2+ (v-u)^2}\dd v, \qq \qq (u\in \R)$$
can be  similarly estimated. We further could not find any reference in the literature. 
\vskip 3 pt

It also appears that formula \eqref{th.LW}-(i), which  exactly means 
\begin{eqnarray}\label{LW-half1} 
\frac1{\pi} \int_{-\infty}^{+\infty} \big|\zeta(\frac1{2} +it)\big|^2{n\over n^2+t^2}\, \dd t&=&   \log n + C_1 + o(1),  
\end{eqnarray}
 seems not be  easily obtained by using complex integration.\end{remark}

\vskip 3 pt  The Stepanov space $\mathcal S^2$ is one instance of amalgam. We recall that the weighted amalgam $\ell^q(L^p,w)$ consists of functions $f$ on $\R$ such that
\begin{equation}\|f\|_{p,w,q}= \Big\{ \,\sum_{n\in \Z} \Big(\int_{n}^{n+1} |f(x)|^p w(x) \dd x\Big)^{\frac{q}{p}}\,\Big\}^{\frac{1}{q}}<\infty,
\end{equation}
where $1<p,q<\infty$ and $w$ is a weight function. We end   with the following 
 question.

\begin{problem}\label{prob3}  In which weighted amalgams $\ell^q(L^2,w)$ lies the  Riemann zeta function?
\end{problem}

\vskip 5 pt
  {\it Final Note.}  We recently solved Problem \ref{prob2}. This is the object of a separate writing \cite{W2b}.

\bigskip
\bigskip
\noindent {\bf Acknowledgments} 
\vskip 9 pt  I thank  Christophe Cuny for   readings and   comments on first attempts. I also thank     Benjamin Enriquez for friendly discussions
 around Problem \ref{prob2}.
 I further thank Antanas Laurin\v cikas and Kohji Matsumoto for comments on  this problem.   Finally, I thank J\"orn Steuding for exchanges around Remark \ref{Parseval-zeta-remark}.


\begin{thebibliography}{99}
 \bibitem{A}T. M. Apostol, \emph{Introduction to analytic number theory}. Springer-Verlag, New York-Heidelberg, (1976), Undergraduate Texts in Mathematics.
  \bibitem{B}
 N. K. Bary,  (1964) 
       \emph{A Treatise on Trigonometric Series} Vol I\&II, translation by M. F. Mullins, Pergamon Press, Oxford-London-Edimburgh-New York-Paris-Frankfurt. 
   \bibitem{BD} P. Bateman  and H. G. Diamond,  (2004) \emph{Analytic number theory--An introductory course}, World Scientific Publishing, Singapore.
    \bibitem{Ca}  L. Carlitz, (1960) Some finite summation formulas of arithmetical character. II, \emph{Acta Math. Acad. Sci. Hung.} {\bf 11}, 15--22.
   \bibitem{CP}  P. Codec\`a and A. Perelli, (1988) On the Uniform Distribution  (mod 1) of the Farey  Fractions and $\ell^p$ Spaces, \emph{Math. Ann.} {\bf 279},   413--422.
\bibitem {C} M. W. Coffey,  (2011) Evaluation of some second moment and other integrals for the Riemann, Hurwitz, and Lerch zeta functions,   arXiv:1101.5722v1.
  \bibitem{CW} C. Cuny and M. Weber, (2018) Pointwise convergence of   almost periodic Fourier series and associated series of dilates, \emph{Pacific J. Math.} {\bf 292}, No. 1, 81--101.
  \bibitem{D}
 H. B.    Dwight,  (1961) 
       \emph{Tabels of Integrals and Other Mathematical Data},
   MacMillan, New York. 
  \bibitem{Gr} T. H.  Gr\"onwall,   (1912)
Some asymptotic expressions in the theory of numbers,  Trans. Amer. Math. Soc. {\bf 8},  118--122.  
 \bibitem{HW} G. H. Hardy et E. M. Wright,      \emph{An introduction to the theory of numbers}, fifth edition, Oxford at the Clarendon Press, Oxford, (1979).
\bibitem {I} A. Ivi\'c, \emph{Some identities for the Riemann zeta function II}, Facta Univ. 20, 1-8 (2005); arXiv:0506214v2 (2005).
\bibitem{LW} M.  Lifshits and  M. Weber (2009)  Sampling the Lindel\" of conjecture with the Cauchy random walk,   \emph{Proc. London Math. Soc.}
{\bf 98} (3),    241--270. \bibitem{L}  J. E. Littlewood, (1912)  Quelques cons\'equences de l'hypoth\`ese que la fonction $\zeta(s)$ de Riemann n'a pas de z\'eros dans le demi-plan $\Re s >\frac12$,   \emph{Comptes Rendus Acad. Sci. Paris}, {\bf154},  263--266.
\bibitem{Ma} K. Matsumoto, (2000) Recent Developments in the Mean Square Theory of the Riemann Zeta and Other Zeta-Functions, \emph{Number Theory}, 241--246, \emph{Trends Math.}, Birkh\"auser, Basel. 
\bibitem{M1}  M.  Mikol\'as, (1949) Farey series and their connection with the prime problem
{(I)}, \emph{Acta Sci. Math. Szeged} {\bf 13},   93--117.
 \bibitem{M2}  M.  Mikol\'as, (1949) {\rm   Farey series and their connection with the prime
problem (II)}, \emph{Acta Sci. Math. Szeged,} {\bf 14},   5--21.
 
\bibitem{MV} H. Montgomery and R. Vaughan, \emph{Multiplicative number theory: I. Classical Theory}, Cambridge Studies in Advanced Math. {\bf 97}, (2006), Cambridge UK.
 \bibitem{Pi}  J. Pintz,    (1983) Oscillatory properties of $M(x)=\sum_{n\le x} \m(n)$, III, \emph{Acta Arith} {\bf 43}, 105--113. 
 \bibitem{Te} G. Tenenbaum, (2008) \emph{Introduction \`a la th\'eorie analytique et probabiliste des nombres}, Coll. \'Echelles Ed. Belin  
 Paris.
  \bibitem{Ti}  E. C. Titchmarsh,  (1986)      \emph{ The theory of
the Riemann-Zeta function}, Second Edition, revised by Heath-Brown D. R., Oxford
Science Publications.  
\bibitem{W2}  M. Weber,   (2019) On infinite M\"obius inversion, \emph{Preprint}.
\bibitem{W2a}  M. Weber,   (2019) Sampling the Riemann zeta function with the Cauchy random walk on the vertical line $\s+it$, $\s>\frac12$, \emph{In progress}.
\bibitem{W2b}  M. Weber,   (2019)  The Stepanov norm of the Riemann zeta function, 
\emph{being drafted}.
\bibitem{W3}  M. Weber,   (2016) A remark on 
Jordan double sums and Dirichlet polynomials, \emph{unpublished notes}.
\bibitem{W1}  M. Weber,   (2009)  \emph{Dynamical Systems and Processes}, European Mathematical
Society Publishing House, IRMA Lectures in Mathematics and Theoretical Physics {\bf 14}, (2009). xiii+759p.    
\bibitem{Wi} A. Wintner, (1957)
  {\rm  Fourier constants and equidistant Riemann sums}, \emph{J. Math. Pures Appl.}
{\bf 36},    251--261.
\bibitem{Y}  M. Yoshimoto  (2004) Abelian Theorems, Farey Series and the Riemann Hypothesis     \emph{The Ramanujan Journal} {\bf 8},  131--145. 
   \bibitem{Z} A. Zygmund,  (2002)  \emph{Trigonometric series}, Third Ed. Vol. {\bf  1}\&{\bf  2} combined,
Cambridge Math. Library, Cambridge Univ. Press.  
   
   
   
   
   \end{thebibliography}
\end{document}